\documentclass[10pt]{article}
\usepackage{latexsym}
\usepackage{amsthm}
\usepackage{amsfonts}
\usepackage{amssymb}
\usepackage{amsmath} 
\usepackage{graphicx}
\usepackage{url} 
\usepackage{color}

\usepackage[top=30truemm,bottom=30truemm,left=30truemm,right=30truemm]{geometry}

\theoremstyle{plain}
  \newtheorem{theorem}{Theorem}[section] 
  
  \newtheorem{proposition}{Proposition}[section]
  \newtheorem{lemma}{Lemma}[section]
  \newtheorem{corollary}{Corollary}[section]

\theoremstyle{remark}
  \newtheorem{remark}{Remark}[section]

\theoremstyle{definition}
  \newtheorem{definition}{Definition}[section]

\include{psfig}

\makeatletter

\@addtoreset{equation}{section}
\makeatother

\include{psfig}

\begin{document}

\markboth{Hideo Takaoka}{Nonlinear Schr\"odinger equations}

\title{Energy transfer model and large periodic boundary value problem for the quintic nonlinear Schr\"{o}dinger equations}
\author{Hideo Takaoka\thanks{This work was supported by JSPS KAKENHI Grant Number 18H01129.}\\
Department of Mathematics, Kobe University\\
Kobe, 657-8501, Japan\\
takaoka@math.kobe-u.ac.jp}

\date{\empty}

\maketitle

\begin{abstract}
We study a dynamics and energy exchanges between a linear oscillator and a nonlinear interaction state for the one dimensional, quintic nonlinear Schr\"odinger equation. 
Gr\'ebert and Thomann \cite{gt1} proved that there exist solutions with initial data built on four Fourier modes, that confirms the conservative exchange of wave energy.
Captured multi resonance in multiple Fourier modes, we simulate a similar energy exchange in long-period waves. 
\end{abstract}

{\it $2010$ Mathematics Subject Classification Numbers.}
35Q55, 42B37.


\section{Introduction}\label{intro}
\indent

In this paper, we consider the defocusing quintic nonlinear Schr\"odinger equation
\begin{eqnarray}\label{nls}
i\partial_t u+\partial_{x}^2u=|u|^4u,\quad t\in\mathbb{R},~x\in\mathbb{T}_L=[-L/2,L/2], 
\end{eqnarray}
where $L>0,~u=u(t,x)\,:\,\mathbb{R}\times\mathbb{T}_L\to\mathbb{C}$ is a complex-valued function and the spatial domain $\mathbb{T}_L$ is taken to be a torus of length $L$, i.e., we assume the periodic boundary condition.
In the case when $L=2\pi$, we denote $\mathbb{T}=\mathbb{T}_{2\pi}$ as usual.

Our aim of this paper is to consider the long periodic solutions ($L\gg 1$) to (\ref{nls}), while there are exchanges of resonant energy at particular frequencies. 
The sign of nonlinearity ($+1$ for the defocusing case and $-1$ for the focusing case) will not play the central role in the present discussion.
For simplicity, we focus on the situation of the defocusing case.

The equation (\ref{nls}) satisfies the mass $M[u]$ and energy $E[u]$ conservations laws;
\begin{eqnarray}\label{con:mass}
M[u](t)=\int_{\mathbb{T}_L}|u(t,x)|^2\,dx,
\end{eqnarray}
\begin{eqnarray}\label{con:energy}
E[u](t)=\int_{\mathbb{T}_L}\left(\frac12|\partial_xu(t,x)|^2+\frac{1}{6}|u(t,x)|^6\right)\,dx,
\end{eqnarray}
which impose the constraints on a dynamics of mass density of solutions.

We briefly recall known results concerning the Cauchy problem for the quintic NLS.
In the non-periodic scenario, i.e., $x\in\mathbb{R}$, the equation is called mass-critical or $L^2$-critical from the viewpoint of scaling.
Indeed, the one-dimensional quintic nonlinear Schr\"odinger equations with non-periodic boundary condition is left invariant by the scaling
$$
u\mapsto u_{\lambda};\quad u(t,x)\mapsto u_{\lambda}(t,x)=\lambda^{1/2}u(\lambda^2 t,\lambda x),\qquad \lambda>0,
$$
which preserves the homogeneous Sobolev norm $\dot{H}^s(\mathbb{R})$ with $s=0$. 
On $\mathbb{R}$-case, the local well-posedness was proved by Cazenave and Weissler \cite{cw1} for data in $L^2$ (see also \cite{gv1} and \cite{ts1}).
Notice that in \cite{cw1}, the existence time of solution depends on the position of data and not only on its size.
One can also prove the global well-posedness in $L^2$ provided that the initial data in $L^2$ is sufficiently small.
Concerning the local well-posedness theory in fractional Sobolev spaces, we refer to the paper \cite{gv1}.
The global well-posedness and scattering in the threshold space $L^2$ was obtained by Dodson \cite{d2}.
More precisely, it is shown that for all $u_0\in L^2(\mathbb{R})$, there exist a unique time global solution to (\ref{nls}) and  $u_{\pm}\in L^2(\mathbb{R})$ satisfying that
$$
\|e^{-it\partial_x^2}u(t)-u_{\pm}\|_{L^2}\to 0
$$
as $t\to\pm\infty$.
It is possible to consider the global well-posedness and scattering for the nonlinear Schr\"odinger equations with energy sub-critical nonlinearities.
This is established in \cite{d1} for initial data below the energy norm.

We now turn to the case of periodic boundary conditions.
The local well-posedness was proved by Bourgain \cite{bo1} for data in $H^s(\mathbb{T})$, with $s>0$.
This combined with the $H^1$-energy conservation law (an a priori estimate for solutions) leads to global well-posedness in $H^1(\mathbb{T})$. 
Similar results hold for the equation in the $L$-periodic boundary condition case $\mathbb{T}_L$ for any $L>0$, without having to change the proof.  

\begin{definition}
For the function $\phi:\mathbb{T}_L\to\mathbb{C}$, we define the Fourier transform 
$$
\widehat{\phi}(\xi)=\int_{-L/2}^{L/2}e^{-ix\xi}\phi(x)\,dx,\quad \xi\in 2\pi\mathbb{Z}/L.
$$
Then we have the representation $\phi(x)=\int e^{ix\xi}\widehat{\phi}(\xi)\,(d\xi)_L$ by
$$
\int e^{ix\xi}\widehat{\phi}(\xi)\,(d\xi)_L=\frac{1}{L}\sum_{\xi\in 2\pi\mathbb{Z}/L}e^{ix\xi}\widehat{\phi}(\xi).
$$
The Sobolev $H^s(\mathbb{T}_L)$ norm is given by
$$
\|\phi\|_{H^s(\mathbb{T}_L)}=\left(\frac{1}{L}\sum_{\xi\in 2\pi\mathbb{Z}/L}\langle\xi\rangle^{2s}|\widehat{\phi}(\xi)|^2\right)^{1/2}.
$$
\end{definition}

In this paper, we want to understand the interaction of mass for each frequency $\xi\in 2\pi\mathbb{Z}/L$.
In fact, the $L^2$-norm conservation law (\ref{con:mass}) constrains this object.
For the periodic boundary value problem, i.e., $x\in\mathbb{T}$ (the case of $L=2\pi$) with replacing the nonlinearity $|u|^4u$ by $\nu |u|^4u$ with $\nu>0$
\begin{eqnarray}\label{GT-NLS}
i\partial_tu+\partial_{x}^2u=\nu |u|^4u,\quad (t,x)\in\mathbb{R}\times \mathbb{T},
\end{eqnarray}
Gr\'ebert and Thomann \cite{gt1} examined the dynamics exhibited by the solution of (\ref{GT-NLS}).
More precisely, they proved the following theorem.

\begin{theorem}[Gr\'ebert and Thomann \cite{gt1}]\label{thm:gt}
Let $k\in\mathbb{Z}\backslash\{0\}$ and $n\in\mathbb{Z}$.
${\cal A}$ is a set of the form ${\cal A}=\{a_2,~a_1,~b_2,~b_1\}$ where $a_2=n,~a_1=n+3k,~b_2=n+4k,~b_1=n+k$.
There exist $T>0,~\lambda_0>0$ and a $2T$-periodic function $K_*:\mathbb{R}\mapsto (0,1)$ which satisfies $K_*(0)\le 1/4$ and $K_*(T)\ge 3/4$ so that if $0<\nu<\nu_0$, there exists a solution to (\ref{GT-NLS}) satisfying for all $0\le t\le \nu^{-3/2}$
$$
u(t,x)=\sum_{j\in{\cal A}}u_j(t)e^{ijx}+\nu^{1/4}q_1(t,x)+\nu^{3/2}tq_2(t,x),
$$
with
$$
|u_{a_1}(t)|^2=2|u_{a_2}(t)|^2=K_*(\nu t),
$$
$$
|u_{b_1}(t)|^2=2|u_{b_2}(t)|^2=1-K_*(\nu t),
$$
and where for all $s\in\mathbb{R},~\|q_1(t,\cdot)\|_{H^s(\mathbb{T})}\le C_s$ for all $t\in\mathbb{R}_+$, and $\|q_2(t,\cdot)\|_{H^s(\mathbb{T})}\le C_s$ for all $0\le t\le \nu^{-3/2}$.
\end{theorem}

From Theorem \ref{thm:gt}, we obtained that there exist solutions with initial data built on four Fourier modes,  that involves periodic energy exchanges between the modes initially excited.
The proof of this result is mostly by calculations the resonant normal form of the Hamiltonian of (\ref{GT-NLS}) up to order ten.

\begin{remark}
Another interesting result is the two-dimensional cubic nonlinear Schr\"odinger equation.
In \cite{ckstt2}, Colliander, Keel, Staffilani, Takaoka and Tao showed the weak turbulence property for the 2D defocusing cubic nonlinear Schr\"odinger equation:
\begin{eqnarray}\label{2dnls}
i\partial_tu+\Delta u=|u|^2u,\qquad (t,x)\in\mathbb{R}\times\mathbb{T}^2.
\end{eqnarray}
More precisely, for any $s>1,~K\gg 1\gg \varepsilon>0$, there exists a time $T\gg 1$ such that the initial value problem corresponding to the equation in (\ref{2dnls}) has a global in time solution $u(t)$ satisfying that
$$
\|u(0)\|_{H^s(\mathbb{T}^2)}\le \varepsilon,\quad \|u(T)\|_{H^s(\mathbb{T}^2)}\ge K.
$$
This
exhibits the $H^s$-norm inflation of solutions to the cubic nonlinear Schr\"odinger equations, that admits solutions with transferring wave energy  from low to high Fourier modes.
\end{remark}

In this paper, we proved that there exist solutions of the one-dimensional quintic nonlinear Schr\"odinger equations with initially excited in multi-frequency modes, where the mass is localized and involves conservative energy exchange between the modes initially excited.
Let us now define the wavenumber set consisting of nonlinear resonance interactions in the equation (\ref{nls}).
\begin{definition}[Resonance interaction set]
Let $k\in\mathbb{N}/L$ be fixed.
For $j\in\mathbb{Z}$, we set $\alpha_{1,j},~\alpha_{3,j},~\alpha_{2,j}$ and $\alpha_{4,j}$ as follows: 
$$
\alpha_{1,j}=2\pi\left(3k+\frac{j}{L}\right),~\alpha_{3,j}=2\pi\frac{j}{L},~\alpha_{2,j}=2\pi\left(k+\frac{j}{L}\right),~\alpha_{4,j}=2\pi\left(4k+\frac{j}{L}\right).
$$
With $\alpha_{m,j}$, we set
$$
{\cal R}_{m}=\left\{\alpha_{m,j}\mid j\in\mathbb{Z},~0\le j< L\right\},
$$
for $1\le m\le 4$, and ${\cal R}=\cup_{m=1}^4{\cal R}_{m}$.
\end{definition}

We start by defining the following norms.
\begin{definition}
We define the sets ${\cal N}_l,~{\cal N}_m,~{\cal N}_r,~{\cal N}_h$ to be subsets of the set $2\pi \mathbb{Z}/L$ as follows:
$$
{\cal N}_l=\left\{\xi\in 2\pi \mathbb{Z}/L   \mid  \xi=2\pi (k\eta+\tau+j/L),~(\eta,\tau,j)\in \mathbb{Z}^3,~\tau\in[0,k),~j\in[0,L),~-99\le \eta\le 98\right\},
$$
$$
{\cal N}_m=\left\{\xi\in 2\pi \mathbb{Z}/L   \mid  \xi=2\pi (k\eta+\tau+j/L),~(\tau,j)\in \mathbb{Z}^2,~\tau\in[0,k),~j\in[0,L),~\eta\in\{99,-100\}\right\},
$$
$$
{\cal N}_r=\left\{\xi\in 2\pi \mathbb{Z}/L   \mid  \xi=2\pi (k\widetilde{\eta}+\widetilde{\tau}+j/L),~(\widetilde{\tau},j)\in \mathbb{Z}^2,~-k/2<\widetilde{\tau}\le k/2,~j\in[0,L),~\widetilde{\eta}\in\{0,1,3,4\}\right\},
$$
$$
{\cal N}_h=(2\pi \mathbb{Z}/L)\backslash({\cal N}_l\cup{\cal N}_m).
$$ 

For the Sobolev index $s\in\mathbb{R}$, let $m(\xi)$ be the multiplier function defined on $2\pi\mathbb{Z}/L$ to be
\begin{eqnarray*}
m(\xi)=\left\{
\begin{array}{ll}
\langle\widetilde{\tau}\rangle^{s-1/2}, & \mbox{if $\xi\in {\cal N}_r$},\\
\langle k\rangle^{s-1/2}, & \mbox{if $\xi\in {\cal N}_l$ and $\xi\not\in{\cal N}_r$},\\
\langle k\rangle^{s-1/2}\langle \tau\rangle^{1/2}, & \mbox{if $\xi\in {\cal N}_m$ and $\xi>0$},\\
\langle k\rangle^{s-1/2}\langle k-\tau\rangle^{1/2}, & \mbox{if $\xi\in {\cal N}_m$ and $\xi<0$},\\
\langle \xi\rangle^s, & \mbox{if $\xi\in {\cal N}_h$},
\end{array}
\right.
\end{eqnarray*}
where $\tau$ and $\widetilde{\tau}$ obey the formula
\begin{eqnarray*}
\xi=2\pi \times
\left\{
\begin{array}{ll}
\left(k\widetilde{\eta}+\widetilde{\tau}+j/L\right), & (\widetilde{\tau},j)\in \mathbb{Z}^2,~-k/2<\widetilde{\tau}\le k/2,~j\in[0,L),~\widetilde{\eta}\in\{0,1,3,4\},\\
\left(k\eta+\tau+j/L\right), & \mbox{otherwise}.
\end{array}
\right.
\end{eqnarray*}
For the sequence $u=(u_{\xi})_{\xi\in 2\pi\mathbb{Z}/L}$, define
$$
\|u\|_s=\left(\frac{1}{L}\sum_{\xi\in 2\pi\mathbb{Z}/L}m(\xi)^2|u_{\xi}|^2\right)^{1/2}.
$$
\end{definition}

\begin{remark}
The set ${\cal R}$ is contained in ${\cal N}_r$.
\end{remark}

Our main result is the following theorem.

\begin{theorem}\label{thm-main}
Let $s\in (1,3/2)$, let $\nu>0$ be a small constant and let $L$ be a natural number.
Then there exist a positive number $k_0$, so that for $k_0\le k\in \mathbb{N}/L$, there exist a smooth global solution $u(t)$ to (\ref{nls}) such that for $|t|\ll 1/(L^3\nu^2)$,
\begin{eqnarray}\label{approx}
u(t)=\sum_{m=1}^4u_{{\cal R}_m}(t)+e(t),
\end{eqnarray}
where $u_{{\cal R}_m}$ is whose Fourier transform is supported on the frequency space ${\cal R}_m$ so that
$$
\|u_{{\cal R}_1}(t)\|_{L^2}=2\|u_{{\cal R}_3}(t)\|_{L^2}^2=\nu\left(\frac{1}{2}-\nu^2K(t)\right),
$$
and
$$
\|u_{{\cal R}_2}(t)\|_{L^2}^2=2\|u_{{\cal R}_4}(t)\|_{L^2}^2=\nu\left(\frac{1}{2}+\nu^2K(t)\right),
$$
where
$$
K(t)=\frac{1}{2\nu^2}\sin\arctan\left(\frac{3\nu^2 t}{2L^3}\right).
$$
Moreover, the error term $e(t)$ associated with the expression in (\ref{approx}) satisfies that
$$
\|e(t)\|_s^2\lesssim \frac{\nu}{L}\left(e^{c\nu^2 t/L^3}-1\right)+ \frac{\nu^3}{\langle k\rangle^{5/2-s}}e^{c\nu^2 t}+\frac{\nu^3}{\langle k\rangle^{3/2-s}}\left(e^{c\nu^2 t}-1\right)+\nu\left(e^{c\nu^2t}-1-c\nu^2t\right).
$$
\end{theorem}

\begin{remark}
In Theorem \ref{thm-main}, we pick $L$ to be a natural number for simplicity.
Only a small modification for the proof is required in the general $L>0$.
It is not important to assume the natural number $L$. 
\end{remark}

\begin{remark}\label{remark:n=0}
The equation (\ref{nls}) has a gauge symmetry.
After conducting gauge transformation $u\mapsto e^{i\theta}u$ with $\theta=n/L$, it is easily to shift $\alpha_{m,j}\mapsto \alpha_{m,j}+2\pi n/L~(1\le m\le 4)$ for $n\in\mathbb{Z}$ and prove the corresponding theorem to Theorem \ref{thm-main}. 
\end{remark}


\begin{remark}
Choosing $k$ large enough such as $k^{3/2-s}\gg L^6\nu^2$, one can prove by Theorem \ref{thm-main} that
$$
\|e(T_0)\|_s^2\ll \left(\frac{\nu}{L}\right)^3T_0,
$$
where $T_0=o(1)/(L^3\nu^2)$.
Then for small data such that $\nu\ll L^{-3/2}$, we have that $T_0\gg 1$ and a solution $u(t)$ to (\ref{nls}) satisfying
$$
\|u(T_0)\|_{H^s}^2-\|u(-T_0)\|_{H^s}^2 \gtrsim  \langle k\rangle^{2s}\left(-3^{2s}+1^{2s}+\frac{4^{2s}}{2}\right)\nu^3K(T_0),
$$
for $s\in(1,3/2)$, where the right-hand side is positive.
\end{remark}

\begin{remark}
We can expect similar results to hold for the following more general nonlinearities with essentially the same proof:
$$
i\partial_t u+\partial_{x}^2u=\sum_{j=1}^Ja_{j}|u|^{2j}u
$$
where $a_j\in\mathbb{R}$.
\end{remark}

The proof of Theorem \ref{thm-main} relies on obtaining the dynamics in a toy model (finite dimensional approximation) of nonlinear Schr\"odinger equations along with error estimates between finite dimensional model and the full infinite dimensional model.

In Section \ref{sec:notation}, we present some notation.
Section \ref{sec:rnls} describes the reductions of the equation (\ref{nls}) as an infinite system of ODE's.
In Section \ref{sec:toymodel}, we construct the appropriate Toy model equation associated with the finite dimensional ODE system of reduced NLS equation, and solve it.
In Section \ref{sec:dyn}, we give qualitative estimates for solutions to the finite dimensional ODE system.
In Section \ref{sec:approx}, we prove that the ODE system derived in Section \ref{sec:toymodel} approximates the dynamics of the quintic nonlinear Schr\"odinger equation (\ref{nls}).
Theorem \ref{thm-main} is established in Section \ref{sec:thm}. 

\section{Notation}\label{sec:notation}
\indent

Let us introduce some notation.
We prefer to use the notation $\langle\cdot\rangle=(1+|\cdot|^2)^{1/2}$.
The over dot $\dot{a}(t)$ denotes the derivative of $a(t)$ with respect to time $t$. 

We use $c,~C$ to denote various constants.
We use $A\lesssim B$ to denote $A\le CB$ for some constant $C>0$.
Similarly, we write $A\ll B$ to mean $A\le c B$ for some small constant $c>0$.


For an odd natural number $n$ and complex-valued functions $f_1,~f_2,\ldots,~f_n$ defined on the set $2\pi \mathbb{Z}/L$, we write
the discrete convolution (convolution sum) $[f_1*f_2*\ldots *f_n](\xi)$ as
$$
[f_1*f_2*\ldots *f_n](\xi)=\frac{1}{L^{n-1}}\sum^{*}\prod_{j=1}^nf_1(\xi_j),
$$
where the superscript notation $*$ of $\sum^*$ indicates a sum running over the hyperplane set $\xi_1-\xi_2+\xi_3-\ldots +\xi_n=\xi$ with summation index.

\section{Reductions of the equation (\ref{nls}) as an infinite system of ODEs}\label{sec:rnls}
\indent

In this section, we consider the smooth solution to (\ref{nls}).
Namely, suppose that $u$ is a smooth global in time solution to (\ref{nls}).
Let us start with the ansatz
$$
u(t,x)=\int a_{\xi}(t)e^{ix\xi-it\xi^2}\,(d\xi)_{L}.
$$
In what follows, we shall omit the time variable $t$ of $a_{\xi}(t)$ in abbreviated form without confusion. 
With this transform, the equation (\ref{nls}) becomes
\begin{eqnarray}\label{rfnls}
i\dot{a}_{\xi}=\frac{1}{L^4}\sum^*a_{\xi_1}\overline{a_{\xi_2}}a_{\xi_3}\overline{a_{\xi_4}}a_{\xi_5}e^{-it\phi(\xi_1,\xi_2,\xi_3,\xi_4,\xi_5,\xi)},
\end{eqnarray}
where the factor $\phi(\xi_1,\xi_2,\xi_3,\xi_4,\xi_5,\xi)$ means the oscillation frequency such that
$$
\phi(\xi_1,\xi_2,\xi_3,\xi_4,\xi_5,\xi_6)=\xi_1^2-\xi_2^2+\xi_3^2-\xi_4^2+\xi_5^2-\xi_6^2.
$$

In the next section, we seek the quintic resonant structure in the last term on the right-hand side of (\ref{rfnls}). 

\section{Toy model}\label{sec:toymodel}
\noindent

In this section, we construct a finite dimensional model, in some sense, as the corresponding approximation model of (\ref{rfnls}).
Consequently, it will approximate an exact solution to (\ref{rfnls}).

Let us begin by recalling some earlier results by Gr\'ebert-Thomann \cite{gt1}.
The first one is some arithmetical result.
\begin{lemma}[Lemma 2.1 in \cite{gt1}]\label{r1}
Assume that $(j_1,j_2,j_3,\ell_1,\ell_2,\ell_3)\in (\mathbb{Z}/L)^6$ satisfy
\begin{eqnarray}\label{res}
\left\{
\begin{array}{l}
j_1+j_2+j_3=\ell_1+\ell_2+\ell_3,\\
j_1^2+j_2^2+j_3^2=\ell_1^2+\ell_2^2+\ell_3^2,
\end{array}
\right.
\quad
\mbox{and}
\quad
\{j_1,j_2,j_3\}\ne \{\ell_1,\ell_2,\ell_3\}.
\end{eqnarray}
Then $\{j_1,j_2,j_3\}\cap\{\ell_1,\ell_2,\ell_3\}=\emptyset$.
\end{lemma}

The second result is small cardinality of (\ref{res}).
\begin{lemma}[Lemma 2.2 in \cite{gt1}]\label{r2}
Assume that there exist $(j_1,j_2,j_3,\ell_1,\ell_2,\ell_3)\in(\mathbb{Z}/L)^6$ which satisfy (\ref{res}).
Then the cardinal number of the set satisfying (\ref{res}) is greater than or equal to $4$.
\end{lemma}

Taking into the observation in Lemma \ref{r2}, we have that the cardinal of the pair $(j_1,\ell_1,j_2,\ell_2,j_3,\ell_3)$ without the case $\{j_1,j_2,j_3\}\ne \{\ell_1,\ell_2,\ell_3\}$, appearing in the right-hand side of (\ref{res}) is greater than or equal to $4$.
  
Now we define the resonance interaction sets.
\begin{definition}\label{def:resonance}
We shall say that the pair $(\xi_1,\xi_2,\xi_3,\xi_4,\xi_5,\xi_6)$ satisfies the resonance condition, if the following conditions hold:
\begin{itemize}
\item[(i)]
$\xi_1+\xi_3+\xi_5=\xi_2+\xi_4+\xi_6$,
\item[(ii)]
$(m_1,m_2,m_3,m_4)\in \{1,2,3,4\}^4$ satisfy $\{(m_1,m_3),\,(m_2,m_4)\}=\{(1,3),\,(2,4)\}$,
\item[(iii)]
two elements of $\xi_1,\xi_3,\xi_5$ are permitted in the set ${\mathcal R}_{m_1}$, namely that are described by $\alpha_{m_1,j_1}$ and $\alpha_{m_1,j_3}$ with some integers $j_1$ and $j_3$; one element of $\xi_1,\xi_3,\xi_5$ is permitted in the set ${\mathcal R}_{m_3}$, namely that is described by $\alpha_{m_3,j_5}$ with some integer $j_5$,
\item[(iv)]
two elements of $\xi_2,\xi_4,\xi_6$ are permitted in the set ${\mathcal R}_{m_2}$, namely that are described by $\alpha_{m_2,j_2}$ and $\alpha_{m_2,j_4}$ with some integers $j_2$ and $j_4$; one element of $\xi_2,\xi_4,\xi_6$ is permitted in the set ${\mathcal R}_{m_4}$, namely that is described by $\alpha_{m_4,j_6}$ with some integer $j_6$,
\item[(v)]
$j_1,j_3,j_5,j_2,j_4,j_6$ given by (ii) and (iii) above satisfy that $\{j_1,j_3\}=\{j_2,j_4\}$ and $j_5=j_6=(j_1+j_3)/2=(j_2+j_4)/2$.
\end{itemize}
\end{definition}

For the resonant condition in Definition \ref{def:resonance}, we have the following property.
\begin{lemma}
If the pair $(\xi_1,\xi_2,\xi_3,\xi_4,\xi_5,\xi_6)$ satisfies the resonance condition, then $\phi(\xi_1,\xi_2,\xi_3,\xi_4,\xi_5,\xi_6)=0$.
\end{lemma}
\noindent
{\it Proof.}
By symmetry, we may assume that $\xi_m~(1\le m\le 6)$ satisfy
$$
(\xi_1,\xi_2,\xi_3,\xi_4,\xi_5,\xi_6)=(\alpha_{1,j_1},\alpha_{2,j_2},\alpha_{1,j_3},\alpha_{2,j_4},\alpha_{3,j_5},\alpha_{4,j_6}),
$$
where
\begin{eqnarray}\label{res:add}
\{j_1,j_3\}=\{j_2,j_4\},\quad j_5=j_6=\frac{j_1+j_3}{2}=\frac{j_2+j_4}{2},
\end{eqnarray}
which yields that $\phi(j_1,j_2,j_3,j_4,j_5,j_6)=0$.
By means of the identity
\begin{eqnarray}\label{res:GT}
\alpha_{1,0}^2+\alpha_{1,0}^2+\alpha_{3,0}^2=\alpha_{2,0}^2+\alpha_{2,0}^2+\alpha_{4,0}^2,
\end{eqnarray}
we obtain
$$
\frac{1}{(2\pi)^2}\phi(\xi_1,\xi_2,\xi_3,\xi_4,\xi_5,\xi_6)=\frac{2k}{L}\left(3(j_1+j_3)-(j_2+j_4)-4j_6\right)+\frac{1}{L^2}\phi(j_1,j_2,j_3,j_4,j_5,j_6),
$$
which equals to zero by (\ref{res:add}).
\qed

\begin{remark}
If we stipulate the $2\pi$-periodic setting to (\ref{nls}), namely $L=2\pi$, such a resonance condition was already known in \cite{gt1}.
Indeed, by taking $\{\xi_1,\xi_3,\xi_5\}=\{\alpha_{1,0},\alpha_{1,0},\alpha_{3,0}\}$ and $\{\xi_2,\xi_4,\xi_6\}=\{\alpha_{2,0},\alpha_{2,0},\alpha_{4,0}\}$ ($j$-factor of $\alpha_{m,j}$ are zero), we obtain
$$
\phi(\xi_1,\xi_2,\xi_3,\xi_4,\xi_5,\xi_6)=\alpha_{1,0}^2+\alpha_{1,0}^2+\alpha_{3,0}^2-\alpha_{2,0}^2-\alpha_{2,0}^2-\alpha_{4,0}^2=0,
$$
which was used in (\ref{res:GT}).
We note that the condition (v) in Definition \ref{def:resonance} implies $\{j_1,j_3,j_5\}=\{j_2,j_4,j_6\}$, which case appears in the second condition in (\ref{res}).
\end{remark}

By removing non-resonance term in the term on the right-hand side of (\ref{rfnls}), we propose to a finite dimensional system of ODE, which  we call the resonant system corresponding to (\ref{rfnls}).
Consider the initial value problem for the following resonant truncation of (\ref{rfnls});
for $\xi=\alpha_{m,j}\in{\mathcal R}_m,~1\le m\le 4$,
\begin{eqnarray}\label{RFNLS}
i\dot{r}_{\xi}& = &
\frac{1}{L^4}\sum^*_{\text{res}(\xi)}r_{\xi_1}\overline{r_{\xi_2}}r_{\xi_3}\overline{r_{\xi_4}}r_{\xi_5},
\end{eqnarray}
where we 
denote by
$\text{res}(\xi)$ the set such that the pair $(\xi_1,\xi_2,\xi_3,\xi_4,\xi_5,\xi)$ satisfies the resonance condition.

We then prove that the approximate solutions of (\ref{RFNLS}).

\section{Dynamics of approximate solutions for some time}\label{sec:dyn}
\indent

In this section, we shall study the resonant truncation ODE system in (\ref{RFNLS}).

\subsection{Conserved quantities}
\noindent

By the subscript ``$\text{res}$" to the summation on the hyperplane $\xi_1+\xi_3+\xi_5=\xi_2+\xi_4+\xi_6$, we use the summation formula
$$
\sum_{\text{res}}f(\xi_1,\xi_2,\xi_3,\xi_4,\xi_5,\xi_6)  =  \sum_{\scriptstyle \xi_1+\xi_3+\xi_5=\xi_2+\xi_4+\xi_6 \atop{\scriptstyle (\xi_1,\xi_2,\xi_3,\xi_4,\xi_5,\xi_6)\,\text{satisfies the  resonance condition}}}f(\xi_1,\xi_2,\xi_3,\xi_4,\xi_5,\xi_6),
$$
which is equivalent to
$$
\sum_{\xi_1\in{\cal R}}\sum^*_{\text{res}(\xi_1)}f(\xi_1,\xi_2,\xi_3,\xi_4,\xi_5,\xi_6)
$$
and its symmetry-equivalent formula. 

The resonant truncation ODE system in (\ref{RFNLS}) retains conserved quantities as follows.
\begin{lemma}\label{con}
Let $\{r_{\xi}(t)\}$ be a global in time solution to (\ref{RFNLS}).
Then we have the relation:
\begin{eqnarray}\label{L^2}
\frac{d}{dt}\sum_{\xi\in{\mathcal R}}|r_{\xi}(t)|^2=0.
\end{eqnarray}
\end{lemma}
\begin{remark}
The $\ell^2$-norm $(\sum_{\xi\in{\mathcal R}}|r_{\xi}(t)|^2)^{1/2}$ is a conserved quantity for  (\ref{RFNLS}).
\end{remark}
\noindent
{\it Proof of Lemma \ref{con}.}
It will be convenience to raise the frequency representation $\xi$ in (\ref{L^2}) to $\xi_6$. 
Namely it suffices to show that
\begin{eqnarray*}
\frac{d}{dt}\sum_{\xi_6\in{\mathcal R}}|r_{\xi_6}(t)|^2=0.
\end{eqnarray*}
Multiplying $\overline{r_{\xi_6}}$ to (\ref{RFNLS}) and taking the imaginary part, we have
\begin{eqnarray*}
\Im(i\dot{r}_{\xi_6}\overline{r_{\xi_6}})
=\frac{1}{L^4}\Im\sum^*_{\text{res}(\xi_6)}r_{\xi_1}\overline{r_{\xi_2}}r_{\xi_3}\overline{r_{\xi_4}}r_{\xi_5}\overline{r_{\xi_6}},
\end{eqnarray*}
where $\xi_6\in {\cal R}$.
The term on the left-hand side of the above equation will be
$$
-\frac{1}{2}\frac{d}{dt}|r_{\xi_6}|^2.
$$
Then after the summation over $\xi_6\in{\mathcal R}$, we arrive at the following:
\begin{eqnarray*}
-\frac{1}{2}\frac{d}{dt}\sum_{\xi_6\in {\mathcal R}}|r_{\xi_6}(t)|^2
=\frac{1}{2i L^4}\sum_{\xi_6\in{\mathcal R}}\sum^*_{\text{res}(\xi_6)}\left(r_{\xi_1}\overline{r_{\xi_2}}r_{\xi_3}\overline{r_{\xi_4}}r_{\xi_5}\overline{r_{\xi_6}}-
\overline{r_{\xi_1}}r_{\xi_2}\overline{r_{\xi_3}}r_{\xi_4}\overline{r_{\xi_5}}r_{\xi_6}\right),
\end{eqnarray*}
which is zero, since by symmetrization schemes.

\qed

We now turn our attention to the dynamical structures of the resonant truncation ODE system in (\ref{RFNLS}).
In a similar way to \cite{ckstt2,gt1,ta1}, we may rewrite (\ref{RFNLS}) to the equation of motion in the Hamiltonian symplectic coordinates.

If we set 
$$
r_{\xi}(t)=\sqrt{I_{\xi}(t)}e^{i\theta_{\xi}(t)},\quad I_{\xi}>0,~\theta_{\xi}\in\mathbb{R},
$$
we can write (\ref{RFNLS}) as a system.
Actually, inserting this back into (\ref{RFNLS}), we see that for $\xi=\alpha_{m,j}\in{\cal R}$
\begin{eqnarray*}
\frac{i}{2}\frac{\dot{I_{\xi}}}{\sqrt{I_{\xi}}}-\sqrt{I_{\xi}}\dot{\theta_{\xi}}
=\frac{1}{L^4}\sum^*_{\text{res$(\xi)$}}\sqrt{I_{\xi_1}I_{\xi_2}I_{\xi_3}I_{\xi_4}I_{\xi_5}}e^{i(\theta_{\xi_1}- \theta_{\xi_2}+\theta_{\xi_3}-\theta_{\xi_4}+\theta_{\xi_5}-\theta_{\xi})}.
\end{eqnarray*}
Taking the real part, we get that for $\xi=\alpha_{m,j}\in{\cal R}$
\begin{eqnarray}
-\dot{\theta_{\xi }} 
=\frac{1}{L^4}\sum^*_{\text{res$(\xi)$}}\sqrt{\frac{I_{\xi_1}I_{\xi_2}I_{\xi_3}I_{\xi_4}I_{\xi_5}}{I_{\xi}}}\cos\left(\theta_{\xi_1}- \theta_{\xi_2}+\theta_{\xi_3}-\theta_{\xi_4}+\theta_{\xi_5}-\theta_{\xi}\right).\label{theta}
\end{eqnarray}
We take the imaginary part to get 
\begin{eqnarray}\label{I}
\frac{1}{2}\dot{I_{\xi}} 
=\frac{1}{L^4}\sum^*_{\text{res$(\xi)$}}\sqrt{I_{\xi_1}I_{\xi_2}I_{\xi_3}I_{\xi_4}I_{\xi_5}I_{\xi}}\sin\left(\theta_{\xi_1}- \theta_{\xi_2}+\theta_{\xi_3}-\theta_{\xi_4}+\theta_{\xi_5}-\theta_{\xi}\right)
\end{eqnarray}
for $\xi=\alpha_{m,j}\in{\cal R}$.
 
\begin{remark}
The ODE system (\ref{theta})-(\ref{I}) enjoys the symmetry $(\theta_{\alpha_{m,j}},I_{\alpha_{m,j}})\to (\theta_{\alpha_{m,L-1-j}},I_{\alpha_{m,L-1-j}})$.
If the data satisfy
$$
I_{\alpha_{m,j}}(0)=I_{\alpha_{m,L-1-j}}(0),\quad \theta_{\alpha_{m,j}}(0)=\theta_{\alpha_{m,L-1-j}}(0),
$$
then the solutions to (\ref{theta})-(\ref{I}) ensure that 
$$
I_{\alpha_{m,j}}(t)=I_{\alpha_{m,L-1-j}}(t),\quad \theta_{\alpha_{m,j}}(t)=\theta_{\alpha_{m,L-1-j}}(t).
$$
\end{remark}

The non-degeneracy of solutions $I_{\alpha_{m,j}}$ to (\ref{theta})-(\ref{I}) is needed for carrying out calculations in (\ref{theta})-(\ref{I}). 
We will provide specific solutions satisfying such a condition in the following lemma.  

\begin{lemma}\label{non-deg}
Given a small constant $\nu>0$, let $\theta_{\alpha_{m,j}}(0),~I_{\alpha_{m,j}}(0)~(1\le m\le 4,~0\le j<L)$ be  an initial datum satisfying $I_{\xi}(0)\sim \nu$.
Then the initial value problem (\ref{theta})-(\ref{I}) admits a unique classical solution $(\theta_{\alpha_{m,j}}(t),~I_{\alpha_{m,j}}(t))$ for $|t|\le c\nu^{-2}L^3
$, where $c>0$ is a small constant $c>0$, such that 
\begin{eqnarray}
\max_{\xi\in \cal R}\left(\nu|\theta_{\xi}(t)
-\theta_{\xi}(0)|+|I_{\xi}(t)-I_{\xi}(0)|\right)\ll \nu.
\end{eqnarray}
\end{lemma}
\noindent
{\it Proof.}
Let us note that because of the time reflection invariance, it suffices to consider non-negative time.
A bootstrap (continuity) argument allows us to consider the set
$$
T=\sup\left\{t\ge 0\mid \max_{\xi\in{\cal R}}\left(\nu|\theta_{\xi}(t)-\theta_{\xi}(0)|+|I_{\xi}(t)-I_{\xi}(0)|\right)\ll \nu\right\}.
$$
Then it suffices to show $T\ge  c\nu^{-2}L^3$.
Clearly, for $0\le t\le T$, the equations (\ref{theta})-(\ref{I}) may be replaced by
\begin{eqnarray*}
-\dot{\theta}_{\xi} =
O\left(\frac{\nu^2}{L^3}\right),\quad 
\frac{1}{2}\dot{I_{\xi}} =
O\left(\frac{\nu^3}{L^3}\right).
\end{eqnarray*}
Suppose $T\ll \nu^{-2}
L^3$.
Then the continuity of the flow implies
$$
\nu\lesssim \max_{\xi\in{\cal R}}\left(|I_{\xi}(T)-I_{\xi}(0)|+\nu|\theta_{\xi}(T)-\theta(0)|\right)  \lesssim 
\frac{\nu^3}{L^3}T
\ll  \nu,
$$
which gives a contradiction.
Then $T\ge c \nu^{-2}
L^3$ for some constant $c>0$.
\qed

The solutions to (\ref{I}) have the following conserved quantities.
\begin{lemma}\label{lem:con-2}
Let $I_{\xi}>0~(\xi\in {\cal R})$ be the solutions of (\ref{theta})-(\ref{I}).
Then
\begin{eqnarray}\label{con-1}
\frac{d}{dt}
\left(I_{\alpha_{3,j}}(t)+I_{\alpha_{4,j}}(t)\right)=0,\quad
\frac{d}{dt}
\left(I_{\alpha_{1,j}}(t)+I_{\alpha_{2,j}}(t)\right)=0,
\end{eqnarray}
\begin{eqnarray}\label{con-2}
\frac{d}{dt}
\left(I_{\alpha_{1,j}}(t)-2I_{\alpha_{3,j}}(t)\right)=0,\quad 
\frac{d}{dt}
\left(I_{\alpha_{2,j}}(t)-2I_{\alpha_{4,j}}(t)\right)=0,
\end{eqnarray}
\end{lemma}
\noindent
{\it Proof.}
First prove (\ref{con-1}). 
We recall the equations to $I_{\alpha_{3,j}}(t)$ and $I_{\alpha_{4,j}}(t)$ 
 as follows:
\begin{eqnarray*}
& &\frac{2}{L^4}\sum^*_{\text{res}(\alpha_{3,j})}\sqrt{I_{\xi_1}I_{\xi_2}I_{\xi_3}I_{\xi_4}I_{\xi_5}I_{\alpha_{3,j}}}\sin\left(\theta_{\xi_1}- \theta_{\xi_2}+\theta_{\xi_3}-\theta_{\xi_4}+\theta_{\xi_5}-\theta_{\alpha_{3,j}}\right)\\
& = & \frac{6}{L^4}\left(2\sum_{\scriptstyle j=\frac{j_2+j_4}{2} \atop{\scriptstyle j\ne j_2}}+\sum_{j=j_2=j_4}\right)\sqrt{I_{\alpha_{2,j_2}}I_{\alpha_{2,j_4}}I_{\alpha_{4,j}}I_{\alpha_{1,j_2}}I_{\alpha_{1,j_4}}I_{\alpha_{3,j}}}\\
& & \sin\left(\theta_{\alpha_{2,j_2}}+\theta_{\alpha_{2,j_4}}+\theta_{\alpha_{4,j}}-\theta_{\alpha_{1,j_2}}-\theta_{\alpha_{1,j_4}}-\theta_{\alpha_{3,j}}\right),
\end{eqnarray*}
and
\begin{eqnarray*}
& & \frac{2}{L^4}\sum^*_{\text{res}(\alpha_{4,j})}\sqrt{I_{\xi_1}I_{\xi_2}I_{\xi_3}I_{\xi_4}I_{\xi_5}I_{\alpha_{4,j}}}\sin\left(\theta_{\xi_1}- \theta_{\xi_2}+\theta_{\xi_3}-\theta_{\xi_4}+\theta_{\xi_5}-\theta_{\alpha_{4,j}}\right)\\
& = & \frac{6}{L^4}\left(2\sum_{\scriptstyle j=\frac{j_2+j_4}{2} \atop{\scriptstyle j\ne j_2}}+\sum_{j=j_2=j_4}\right)\sqrt{I_{\alpha_{1,j_2}}I_{\alpha_{1,j_4}}I_{\alpha_{3,j}}I_{\alpha_{2,j_2}}I_{\alpha_{2,j_4}}I_{\alpha_{4,j}}}\\
& & \sin\left(\theta_{\alpha_{1,j_2}}+\theta_{\alpha_{1,j_4}}+\theta_{\alpha_{3,j}}-\theta_{\alpha_{2,j_2}}-\theta_{\alpha_{2,j_4}}-\theta_{\alpha_{4,j}}\right).
\end{eqnarray*}
Substituting for $I_{\alpha_{3,j}}(t),~I_{\alpha_{4,j}}(t)$ in the right-hand side of (\ref{I}) yields
\begin{eqnarray*}
\frac{1}{2} \frac{d}{dt}\left(I_{\alpha_{3,j}}(t)+I_{\alpha_{4,j}}(t) \right)= 0.
\end{eqnarray*}
 
We next prove the second estimate in (\ref{con-1}).
If the pair $(\xi_1,\xi_2,\xi_3,\xi_4,\xi_5,\alpha_{1,j})$ satisfies the resonant condition, then $\{\xi_2,\,\xi_4\}=\{\alpha_{1,j_2},\alpha_{3,j_6}\}$ and $\{\xi_1,\xi_3,\xi_5\}=\{\alpha_{2,j_2},\alpha_{2,j},\alpha_{4,j_6}\}$ for some $j_2$ and $j_6$ such that $j_6=(j_2+j)/2$.
Along the similar process as above, we have that the term right-hand side in (\ref{I}) for $I_{\alpha_{1,j}}(t)$ is restated as follows:
\begin{eqnarray*}
&  & \frac{2}{L^4}\sum^*_{\text{res}(\alpha_{1,j})}\sqrt{I_{\xi_1}I_{\xi_2}I_{\xi_3}I_{\xi_4}I_{\xi_5}I_{\alpha_{1,j}}}\sin\left(\theta_{\xi_1}- \theta_{\xi_2}+\theta_{\xi_3}-\theta_{\xi_4}+\theta_{\xi_5}-\theta_{\alpha_{1,j}}\right)\\
& = & \frac{12}{L^4}\left(2\sum_{\scriptstyle j_6=\frac{j_2+j}{2} \atop{\scriptstyle j\ne j_2}}+\sum_{j=j_2=j_4}\right)\sqrt{I_{\alpha_{2,j_2}}I_{\alpha_{2,j}}I_{\alpha_{4,j_6}}I_{\alpha_{3,j_6}}I_{\alpha_{1,j_2}}I_{\alpha_{1,j}}}\\
& & \sin\left(\theta_{\alpha_{2,j_2}}+\theta_{\alpha_{2,j}}+\theta_{\alpha_{4,j_6}}-\theta_{\alpha_{3,j_6}}-\theta_{\alpha_{1,j_2}}-\theta_{\alpha_{1,j}}\right).
\end{eqnarray*}
On the other hand, if the pair $(\xi_1,\xi_2,\xi_3,\xi_4,\xi_5,\alpha_{2,j})$ satisfies the resonant condition, then $\{\xi_2,\xi_4\}=\{\alpha_{2,j_2},\alpha_{4,j_6}\}$ and $\{\xi_1,\xi_3,\xi_5\}=\{\alpha_{1,j_2},\alpha_{1,j},\alpha_{3,j_6}\}$ for some $j_2$ and $j_6$ such that $j_6=(j_2+j)/2$, which deduces that the  term right-hand side in (\ref{I}) for $I_{\alpha_{2,j}}(t)$ given by (\ref{I}) is restated as follows:
\begin{eqnarray*}
& & \frac{2}{L^4}\sum^*_{\text{res}(\alpha_{2,j})}\sqrt{I_{\xi_1}I_{\xi_2}I_{\xi_3}I_{\xi_4}I_{\xi_5}I_{\alpha_{2,j}}}\sin\left(\theta_{\xi_1}- \theta_{\xi_2}+\theta_{\xi_3}-\theta_{\xi_4}+\theta_{\xi_5}-\theta_{\alpha_{2,j}}\right)\\
& = & \frac{12}{L^4}\left(2\sum_{\scriptstyle j_6=\frac{j_2+j}{2} \atop{\scriptstyle j\ne j_2}}+\sum_{j=j_2=j_4}\right)\sqrt{I_{\alpha_{1,j_2}}I_{\alpha_{1,j}}I_{\alpha_{3,j_6}}I_{\alpha_{2,j_2}}I_{\alpha_{4,j_6}}I_{\alpha_{2,j}}}\\
& & \sin\left(\theta_{\alpha_{1,j_2}}+\theta_{\alpha_{1,j}}+\theta_{\alpha_{3,j_6}}-\theta_{\alpha_{2,j_2}}-\theta_{\alpha_{4,j_6}}-\theta_{\alpha_{2,j}}\right).
\end{eqnarray*}
This proves the second estimate in (\ref{con-1}).

We can also conclude that the estimates in (\ref{con-2}) hold by means of similarity computation as above.
\qed

Let us proceed to construction of the specific solution to (\ref{theta}) and (\ref{I}).
We define
$$
\Phi_{j_1,j_2,j_3,j_4,j_5,j_6}^{l_1,l_2,l_3,l_4,l_5,l_6}(t)=\theta_{\alpha_{l_1,j_1}}(t)-\theta_{\alpha_{1_2,j_2}}(t)+\theta_{\alpha_{l_3,j_3}}(t)-\theta_{\alpha_{l_4,j_4}}(t)+\theta_{\alpha_{l_5,j_5}}(t)-\theta_{\alpha_{l_6,j_6}}(t).
$$
We establish the following lemma. 
\begin{lemma}\label{appPhi}
Let $\Phi_{j_1,j_1,j_2,j_2,j_3,j_3}^{1,2,1,2,3,4}(0)=\pi/2$ for all $j_1,~j_2,~j_3$ satisfying  $j_3=(j_1+j_2)/2$, and suppose that $I_{\xi}(t)\sim \nu>0$ for $|t|<T$.
Then for $|t|<T$,
\begin{eqnarray}\label{Phi_1}
\Phi_{j_1,j_1,j_2,j_2,j_3,j_3}^{1,2,1,2,3,4}(t)=\Phi_{j_1,j_1,j_2,j_2,j_3,j_3}^{1,2,1,2,3,4}(0)
\end{eqnarray}
for all $j_1,~j_2,~j_3$ satisfying  $j_3=(j_1+j_2)/2$.
\end{lemma}
\noindent
{\it Proof.}
It is a straight forward matter to obtain the result. 
Since by $\theta_{\alpha_{l_1,j_1}}(0)-\theta_{\alpha_{1_2,j_2}}(0)+\theta_{\alpha_{l_3,j_3}}(0)-\theta_{\alpha_{l_4,j_4}}(0)+\theta_{\alpha_{l_5,j_5}}(0)-\theta_{\alpha_{l_6,j_6}}(0)=\pi/2$, it follows that from (\ref{theta}),
$$
\max_{m,j}|\theta_{\alpha_m,j}(t)-\theta_{\alpha_m,j}(0)|\lesssim \frac{\nu^2}{L^3}\int_0^t\max_{m,j}|\theta_{\alpha_m,j}(t')-\theta_{\alpha_m,j}(0)|\,dt'.
$$
We apply the Gronwall inequality to get
$$
\theta_{\alpha_m,j}(t)=\theta_{\alpha_m,j}(0),
$$
which implies $\Phi_{j_1,j_1,j_2,j_2,j_3,j_3}^{1,2,1,2,3,4}(t)=\Phi_{j_1,j_1,j_2,j_2,j_3,j_3}^{1,2,1,2,3,4}(0).$
\qed

\subsection{Averaging property}
\noindent

It is natural to expect that the average of the $L$ sums by
$$
\frac{1}{L}\sum_{j=0}^{L-1}I_{\alpha_{m,j}}(t)=\frac{I_{\alpha_{m,0}}(t)+I_{\alpha_{m,1}}(t)+\ldots+I_{\alpha_{m,L-1}}(t)}{L},\quad 1\le m\le 4
$$
approximates the source of a mass located in frequency space ${\cal R}_m$.
By Lemma \ref{appPhi}, we may suppose $\Phi_{j_1,j_1,j_2,j_2,j_3,j_3}^{1,2,1,2,3,4}(t)=\pi/2$ in (\ref{theta})-(\ref{I}).
In a certain sense that should be taken the average of both the left- and right-hand sides of (\ref{I}) with respect to $0\le j<L$, we reformulate the ODE system (\ref{I}) as the following ODE system:
\begin{eqnarray}\label{RA}
\begin{cases}
\displaystyle \dot{I}_{{\cal R}_1} = -\frac{12}{L^3}\sqrt{I_{{\cal R}_4} I_{{\cal R}_2}^2I_{{\cal R}_3} I_{{\cal R}_1}^2},&\\
\displaystyle \dot{I}_{{\cal R}_2} =  \frac{12}{L^3}\sqrt{I_{{\cal R}_4} I_{{\cal R}_2}^2I_{{\cal R}_3} I_{{\cal R}_1}^2},&\\
\displaystyle \dot{I}_{{\cal R}_3} = -\frac{6}{L^3}\sqrt{I_{{\cal R}_4} I_{{\cal R}_2}^2I_{{\cal R}_3} I_{{\cal R}_1}^2},&\\
\displaystyle \dot{I}_{{\cal R}_4} = \frac{6}{L^3}\sqrt{I_{{\cal R}_4} I_{{\cal R}_2}^2I_{{\cal R}_3} I_{{\cal R}_1}^2}.
\end{cases}
\end{eqnarray}
A similar argument in Lemmas \ref{con} and \ref{lem:con-2}  shows that
\begin{eqnarray}\label{con-R1}
\frac{d}{dt}\sum_{m=1}^4I_{{\cal R}_m}(t)=\frac{d}{dt}\left(I_{{\cal R}_3}(t)+I_{{\cal R}_4}(t)\right)= \frac{d}{dt}\left(I_{{\cal R}_1}(t)+I_{{\cal R}_2}(t)\right)=0.
\end{eqnarray}
It is easily to see that from (\ref{RA})\begin{eqnarray}\label{con-R}
\frac{d}{dt}\left(I_{{\cal R}_1}(t)-2I_{{\cal R}_3}(t)\right)=\frac{d}{dt}\left(I_{{\cal R}_2}(t)-2I_{{\cal R}_4}(t)\right)=0,
\end{eqnarray}
which are also the conserved quantities of the system in (\ref{RA}).


Let us now use an expression of the form
\begin{eqnarray*}
I={}^t(I_{{\cal R}_2},\,I_{{\cal R}_4},\,I_{{\cal R}_1},\,I_{{\cal R}_3}),
\end{eqnarray*}
and define new variables
\begin{eqnarray*}
J={}^t(J_1,\,J_2,\,J_3,\,J_4),
\end{eqnarray*}
where
$$
J_1=\frac12 I_{{\cal R}_2},\quad J_2=-\frac12I_{{\cal R}_2}+I_{{\cal R}_4},\quad J_3=I_{{\cal R}_2}+I_{{\cal R}_1},\quad J_4=\frac12I_{{\cal R}_2}+I_{{\cal R}_3}.
$$
Then the Hamiltonian flow in the action-angle coordinates (\ref{RA}) satisfies
\begin{eqnarray*}
J
=
A 
I
\end{eqnarray*}
where
\begin{eqnarray*}
A=\left(
\begin{array}{cccc}
2 & 1 & -2 & -1 \\
0 & 1 & 0 & 0 \\
0 & 0 & 1 & 0 \\
0 & 0 & 0 & 1 
\end{array}
\right)
.
\end{eqnarray*}
We will consider a very simple choice of initial data.
Given  $\nu>0$, by adapting the result in (\ref{con-R1}) and (\ref{con-R}), we specialize the variables of $I_{{\cal R}_m}$'s to the normalization form
$$
I_{{\cal R}_1}(t)+I_{{\cal R}_2}(t)=\nu,\quad I_{{\cal R}_3}(t)+I_{{\cal R}_4}(t)=\frac{\nu}{2},
$$ 
$$
I_{{\cal R}_1}(t)-2I_{{\cal R}_3}(t)=I_{{\cal R}_2}(t)-2I_{{\cal R}_4}(t)=0.
$$
Employing the auxiliary function $K(t)$ such that $\nu^2 |K(t)|\in[0,1/2)$, we put
\begin{eqnarray}\label{K-fun}
I_{{\cal R}_1}(t)=2I_{{\cal R}_3}(t)=\nu\left(\frac{1}{2}-\nu^2 K(t)\right),\quad I_{{\cal R}_2}(t)=2I_{{\cal R}_4}(t)=\nu\left(\frac{1}{2}+\nu^2 K(t)\right),
\end{eqnarray}
which imply that
$$
J_2=0,\quad J_3=\nu,\quad J_4=\frac{\nu}{2}.
$$
Solutions to (\ref{RA}) with the constrain in (\ref{K-fun}) are proposed to solve the following equations on 
$K(t)$:
\begin{eqnarray}\label{pK}
\displaystyle \dot{K}= \frac{6}{L^3}\left(\frac{1}{2}+\nu^2 K\right)^{3/2}\left(\frac{1}{2}-\nu^2 K\right)^{3/2}.
\end{eqnarray}
We remark that a solution of the ordinary differential equation
$$
\dot{f}(t)=a\left(b-f(t)\right)^{3/2}\left(b+f(t)\right)^{3/2}
$$
is provided by
$$
f(t)=b\sin\arctan\left(ab^2t+c\right)
$$
where $c\in\mathbb{R}$ is a constant.
Therefore, we will choose the function $K$ such that
\begin{eqnarray}\label{choice-K}
K(t)=\frac{1}{2\nu^2}\sin\arctan\left(\frac{3\nu^2 t
}{2L^3}\right).
\end{eqnarray}
%

Dealing with the special choice of the function, we use the perturbation theory for finding an approximate solution to (\ref{I}).
We start from the following calculation:
\begin{eqnarray}\label{ap-1}
\sum_{j=0}^{L-1}\mbox{card}\left\{(j_1,j_2)\in\mathbb{Z}^2\mid 0\le j_1,j_2< L,~j=\frac{j_1+j_2}{2}\right\}=\sum_{j=0}^{L-1}\min\{j+1,L-j\}=\frac{L^2}{4}+O(L),
\end{eqnarray}
and
\begin{eqnarray}
& & \sum_{j_2=0}^{L-1}\mbox{card}\left\{(j,j_1)\in\mathbb{Z}^2\mid 0\le j_1,j< L,~j=\frac{j_1+j_2}{2}\right\}\label{ap-2}\\
& = & \sum_{\scriptstyle  j_2=0 \atop{\scriptstyle j_2:\text{even}}}^{L-1}\mbox{card}\left\{(j,j_1)\in\mathbb{Z}^2\mid 0\le j_1,j< L,~j=\frac{j_1+j_2}{2}\right\}\nonumber\\
&   & +\sum_{\scriptstyle  j_2=0 \atop{\scriptstyle j_2:\text{odd}}}^{L-1}\mbox{card}\left\{(j,j_1)\in\mathbb{Z}^2\mid 0\le j_1,j< L,~j=\frac{j_1+j_2}{2}\right\}\nonumber\\
& = & \sum_{\scriptstyle  j_2=0 \atop{\scriptstyle j_2:\text{even}}}^{L-1}\left(\frac{L}{2}+O(1)\right)+\sum_{\scriptstyle  j_2=0 \atop{\scriptstyle j_2:\text{odd}}}^{L-1}\left(\frac{L}{2}+O(1)\right)\nonumber\\
& = & \frac{L^2}{2}+O(L).\nonumber
\end{eqnarray}

\begin{lemma}\label{lem:pert}
Let $(I_{{\cal R}_m})_{1\le m\le 4}$ be global solutions of (\ref{RA}) 
(those solutions were constructed in (\ref{K-fun})-(\ref{choice-K})).
If at $t=0$, the initial datum $(\theta_{\alpha_{m,j}}(0),I_{\alpha_{m,j}}(0))$ satisfies 
$$
I_{\alpha_{m,j}}(0)=I_{{\cal R}_m}(0),
$$
$$
\theta_{\alpha_{m,j}}(0)=\theta_{\alpha_{m}}
$$
for $1\le m\le 4,~0\le j<L$, and
$$
2\theta_{\alpha_{2}}+\theta_{\alpha_{4}}-2\theta_{\alpha_{1}}-\theta_{\alpha_{3}}=\frac{\pi}{2},
$$
then there exist a constant $c>0$ and solutions $(\theta_{m,j},I_{\alpha_{m,j}})$ to (\ref{theta})-(\ref{I}), so that for $|t|\ll 
\nu^{-2}L^3$,
\begin{eqnarray}\label{app}
\left|\frac{1}{L}\sum_{0\le j<L}I_{\alpha_{m,j}}(t)-I_{{\cal R}_m}(t)\right|\lesssim \frac{\nu}{L} \left(e^{c\frac{\nu^2}{L^3}t}-1\right).
\end{eqnarray}
\end{lemma}
\noindent
{\it Proof.}
It suffices to consider non-negative time.
By Lemma \ref{appPhi}, we suppose $\Phi_{j_1,j_1,j_2,j_2,j_3,j_3}^{1,2,1,2,3,4}(t)=\pi/2$.
For fixed $0\le j<L$, use (\ref{ap-2}) to see that
\begin{eqnarray}
& &  \frac{12}{L^4}\left(2\sum_{\scriptstyle j_6=\frac{j_2+j}{2} \atop{\scriptstyle j\ne j_2}}+\sum_{j=j_2=j_4}\right)\sqrt{I_{\alpha_{2,j_2}}I_{\alpha_{2,j}}I_{\alpha_{4,j_6}}I_{\alpha_{3,j_6}}I_{\alpha_{1,j_2}}I_{\alpha_{1,j}}}
 - \frac{12}{L^3}\sqrt{I_{{\cal R}_4} I_{{\cal R}_2}^2I_{{\cal R}_3} I_{{\cal R}_1}^2}\nonumber\\
& = &\frac{24}{L^4}\sum_{\scriptstyle j_6=\frac{j_2+j}{2} \atop{\scriptstyle j\ne j_2}} \left( \sqrt{I_{\alpha_{2,j_2}}I_{\alpha_{2,j}}I_{\alpha_{4,j_6}}I_{\alpha_{3,j_6}}I_{\alpha_{1,j_2}}I_{\alpha_{1,j}}}
- \sqrt{I_{{\cal R}_4} I_{{\cal R}_2}^2I_{{\cal R}_3} I_{{\cal R}_1}^2}
\right)
+O\left(\frac{\nu^3}{L^4}\right).\label{1-1}
\end{eqnarray}
By virtue of the proof of Lemma \ref{lem:con-2}, 
(\ref{1-1}) implies
\begin{eqnarray}
|I_{\alpha_{1,j}}(t)-I_{{\cal R}_1}(t)| & \lesssim & 
\frac{\nu^2}{L^3}\int_0^t\left(\max_{\scriptstyle 1\le m\le 4\atop{\scriptstyle 0\le j<L}}|I_{\alpha_{m,j}}(t')-I_{{\cal R}_m}(t')|
+\frac{\nu}{L}\right)\,dt'.
\label{es:m=1}
\end{eqnarray}
provided $t\ll L^3/\nu^2$. 
Thanks to the conservation laws obtained in Lemma \ref{lem:con-2}, (\ref{con-R1}) and (\ref{con-R}) along with the initial condition at $t=0$, the estimate in (\ref{es:m=1}) implies that
\begin{eqnarray*}
\max_{1\le m\le 4}|I_{\alpha_{m,j}}(t)-I_{{\cal R}_m}(t)|\lesssim  \frac{\nu^2}{L^3}\int_0^t\left(\max_{\scriptstyle 1\le m\le 4\atop{\scriptstyle 0\le j<L}}|I_{\alpha_{m,j}}(t')-I_{{\cal R}_m}(t')|+\frac{\nu}{L}\right)\,dt'.
\end{eqnarray*}
We apply the Gronwall inequality to get
\begin{eqnarray}\label{es:m=2}
\int_0^t\left(\max_{1\le m\le 4}|I_{\alpha_{m,j}}(t')-I_{{\cal R}_m}(t')|+\frac{\nu}{L}\right)\,dt'\lesssim \frac{L^2}{\nu}\left(e^{c\frac{\nu^2}{L^3}t}-1\right).
\end{eqnarray}
%

We control now (\ref{app}).
Use (\ref{ap-1}) to see that
\begin{eqnarray}
& & \frac{1}{L}\sum_{0\le j<L}\frac{6}{L^4}\left(2\sum_{\scriptstyle j=\frac{j_2+j_4}{2} \atop{\scriptstyle j\ne j_2}}+\sum_{j=j_2=j_4}\right)\sqrt{I_{\alpha_{2,j_2}}I_{\alpha_{2,j_4}}I_{\alpha_{4,j}}I_{\alpha_{1,j_2}}I_{\alpha_{1,j_4}}I_{\alpha_{3,j}}}
-\frac{6}{L^3}\sqrt{I_{{\cal R}_4} I_{{\cal R}_2}^2I_{{\cal R}_3} I_{{\cal R}_1}^2}
\nonumber\\
& =  &  \frac{12}{L^5}\sum_{0\le j<L}\sum_{\scriptstyle j=\frac{j_2+j_4}{2} \atop{\scriptstyle j\ne j_2}}\left(\sqrt{I_{\alpha_{2,j_2}}I_{\alpha_{2,j_4}}I_{\alpha_{4,j}}I_{\alpha_{1,j_2}}I_{\alpha_{1,j_4}}I_{\alpha_{3,j}}}
-\sqrt{I_{{\cal R}_4} I_{{\cal R}_2}^2I_{{\cal R}_3} I_{{\cal R}_1}^2}
\right)+O\left(\frac{\nu^3}{L^4}\right).
\label{3-1}
\end{eqnarray}
Using (\ref{es:m=2}) and (\ref{3-1}), 
we arrive at
\begin{eqnarray*}
\max_{\scriptstyle 1\le m\le 4\atop{\scriptstyle 0\le j<L}}\left|\frac{1}{L}\sum_{0\le j<L}I_{\alpha_{m,j}}(t)-I_{{\cal R}_m}(t)\right| & \lesssim  & \frac{\nu^2}{L^3}\int_0^t\left(\max_{\scriptstyle 1\le m\le 4\atop{\scriptstyle 0\le j<L}}|I_{\alpha_{m,j}}(t')-I_{{\cal R}_m}(t')|+\frac{\nu}{L}\right)\,dt'\\
& \lesssim  & \frac{\nu}{L} \left(e^{c\frac{\nu^2}{L^3}t}-1\right),
\end{eqnarray*}
 which completes the proof.
\qed

By (\ref{K-fun}), (\ref{choice-K}) and Lemma \ref{lem:pert}, we will automatically have the following.
\begin{proposition}
Given a small constant $\nu>0$, let the initial datum the initial datum $(\theta_{\alpha_{m,j}}(0),I_{\alpha_{m,j}}(0))$ satisfy that
$$
I_{\alpha_{m,j}}(0)=I_{{\cal R}_m}(0),
$$
$$
\theta_{\alpha_{m,j}}(0)=\theta_{\alpha_{m}}
$$
for $1\le m\le 4,~0\le j<L$, and
$$
2\theta_{\alpha_{2}}+\theta_{\alpha_{4}}-2\theta_{\alpha_{1}}-\theta_{\alpha_{3}}=\frac{\pi}{2},
$$
where $(I_{{\cal R}_m}(0))_{m=1}^4$ are the same as in (\ref{K-fun}), (\ref{choice-K}).
Then for $|t|\ll L^3\nu^{-2}$,
$$
\frac{1}{L}\sum_{0\le j<L}I_{\alpha_{m,j}}(t)=I_{{\cal R}_m}(t)+O\left(\frac{\nu^3}{L^4}t\right).
$$
\end{proposition}

\section{Approximate estimates}\label{sec:approx}
\indent

In this section, we study the approximation of the infinite dimensional NLS flow in (\ref{rfnls}).
In order to pass from the finite-dimensional flow to the infinite one, we shall evaluate residual terms $a_{\xi}(t)$ of $\xi\not\in {\cal R}$ in the infinite-dimensional ODE (\ref{rfnls}) and ultimately approximate the full system (\ref{RFNLS}) corresponding to (\ref{rfnls}).

Suppose that $(a_{\xi}(t))_{\xi\in 2\pi \mathbb{Z}/L}$ and $(r_{\xi}(t))_{\xi\in{\cal R}}$ are solutions to (\ref{rfnls}) and (\ref{RFNLS}), respectively.
As a matter of convenience, we recurse the sequence $(r_{\xi}(t))_{\xi\in 2\pi \mathbb{Z}/L}$ by placing $r_{\xi}(t)=0$ for $\xi\not\in {\cal R}$.
Also we may extend the formula $\text{res}(\xi)$ for all $\xi\in 2\pi \mathbb{Z}/L$ by replacing $\text{res}(\xi)=\emptyset$ for $\xi\not\in {\cal R}$ to deal with the case when $\xi\not\in {\cal R}$.

The initial datum $(a_{\xi}(0),r_{\xi}(0))_{\xi\in 2\pi \mathbb{Z}/L}$ are given as follows:
\begin{itemize}
\item
if $\xi\in{\cal R}_m$ for some $1\le m\le 4$, then $a_{\xi}(0)=r_{\xi}(0)=I_{{\cal R}_m}(0)e^{i\theta_{{\cal R}_m}}$,\\
\item
if $\xi\not\in\cup_{m=1}^4{\cal R}_m$, then $a_{\xi}(0)=r_{\xi}(0)=0$,
\end{itemize}
where $I_{{\cal R}_m}(t)$ are provided in (\ref{K-fun}) and (\ref{choice-K}), and $\theta_{{\cal R}_m}$ satisfy
$$
2\theta_{{\cal R}_2}+\theta_{{\cal R}_4}-2\theta_{{\cal R}_1}-\theta_{{\cal R}_3}=\frac{\pi}{2}.
$$
It follows that from the mass conservation laws (\ref{con:mass}) and (\ref{L^2}),
\begin{eqnarray}\label{ap}
\frac{1}{L}\sum_{\xi\in 2\pi\mathbb{Z}/L}\left(|a_{\xi}(t)|^2+|r_{\xi}(t)|^2\right) \le c \nu,
\end{eqnarray}
where $c>0$ is independent of $t\in\mathbb{R}$.

We define a subset of $2\pi\mathbb{Z}/L$.
\begin{definition}\label{def:tilde}
Introduce some notions.
For $\xi\in 2\pi\mathbb{Z}/L$, rewrite the form 
\begin{eqnarray}\label{expr:xi}
\xi=2\pi \left(k\eta+\tau+\frac{j}{L}\right),~(\eta,\tau,j)\in \mathbb{Z}^3,~\tau\in[0,k),~j\in[0,L)
\end{eqnarray}
as
\begin{eqnarray}\label{tilde}
\xi=2\pi \left(k\widetilde{\eta}+\widetilde{\tau}+\frac{j}{L}\right),
\end{eqnarray}
where
\begin{eqnarray*}
(\widetilde{\eta},\widetilde{\tau})=
\left\{
\begin{array}{ll}
(\eta,\tau), & \mbox{if $\tau\in[0,k/2]$},\\
(\eta+1,\tau-k), & \mbox{if $\tau\in(k/2,k)$}.
\end{array}
\right.
\end{eqnarray*}
With the notions above, define
\begin{eqnarray*}
A_1& = & \cup_{\eta\in\{0,1,3,4\}}\left\{(\xi_1,\xi_2,\xi_3,\xi_4,\xi_5,\xi_6)\in (2\pi\mathbb{Z}/L)^6 \right. \mid   \xi_m=2\pi(k\widetilde{\eta}_m+\widetilde{\tau}_m+j_m/L, \\
& &\left. (\widetilde{\tau}_m,j_m)\in\mathbb{Z}^2,~\widetilde{\eta}_m=k\eta,~-k/2<\widetilde{\tau}_m\le k/2,~j_m\in[0,L),~1\le m\le 6\right\}.
\end{eqnarray*}
\end{definition}

By removing the resonant part of the Fourier transform from the nonlinear interaction of $a_{\xi}(t)$, we write the error of Fourier mode $\xi$ as
$$
e_{\xi}(t)=a_{\xi}(t)-r_{\xi}(t)
$$
for $\xi\in  2\pi \mathbb{Z}/L$.
One uses (\ref{rfnls}) and (\ref{RFNLS}) to obtain something like
\begin{eqnarray}\label{eq:e}
i\dot{e}_{\xi}(t) =\sum_{m=1}^3 R_{\xi}^m(t),
\end{eqnarray}
where (by omitting the time variable $t$ from the equations) 
\begin{eqnarray*}
R_{\xi}^{1} =\frac{1}{L^4}\sum_{(\xi_1,\xi_2,\xi_3,\xi_4,\xi_5,\xi)\in A_1}^*
a_{\xi_1}\overline{a_{\xi_2}}a_{\xi_3}\overline{a_{\xi_4}}a_{\xi_5}e^{-it\phi(\xi_1,\xi_2,\xi_3,\xi_4,\xi_5,\xi)},
\end{eqnarray*}
\begin{eqnarray*}
R_{\xi}^{2} =  \frac{1}{L^4}\sum_{\text{res}(\xi)}^*\left(a_{\xi_1}\overline{a_{\xi_2}}a_{\xi_3}\overline{a_{\xi_4}}a_{\xi_5}-r_{\xi_1}\overline{r_{\xi_2}}r_{\xi_3}\overline{r_{\xi_4}}r_{\xi_5}\right),
\end{eqnarray*}
\begin{eqnarray*}
R_{\xi}^3 = \frac{1}{L^4}\sum^*_{\scriptstyle  \text{res}(\xi)^c\atop{\scriptstyle (\xi_1,\xi_2,\xi_3,\xi_4,\xi_5,\xi)\in A_1^c}}a_{\xi_1}\overline{a_{\xi_2}}a_{\xi_3}\overline{a_{\xi_4}}a_{\xi_5}e^{-it\phi(\xi_1,\xi_2,\xi_3,\xi_4,\xi_5,\xi)}.
\end{eqnarray*}

\begin{remark}\label{rem:existence}
In the spirit of the standard local well-posedness theory observed in \cite{cw1,ts1} and an ODE technique along with the conservation laws in (\ref{con:mass}), (\ref{con:energy}) and (\ref{L^2}), we easily see that there exists a unique smooth global in time solution to the initial value problem for the corresponding equations to (\ref{rfnls}) and (\ref{RFNLS}), respectively. 
\end{remark}

By a slight abuse of notation above, we define the modified $H^s$-energy for the difference between the solutions to (\ref{rfnls}) and (\ref{RFNLS}) as follows:
\begin{eqnarray}\label{m-e}
\widetilde{E}(t)=\|e(t)\|_s^2+\frac{2}{L^5}\Re \sum_{\xi_6}\sum_{\scriptstyle (\xi_1,\xi_2,\xi_3,\xi_4,\xi_5,\xi_6)\in A_1 \atop{\scriptstyle |\phi(\xi_1,\xi_2,\xi_3,\xi_4,\xi_5,\xi_6)|>\delta}}^*m(\xi_6)^2\frac{a_{\xi_1}\overline{a_{\xi_2}}a_{\xi_3}\overline{a_{\xi_4}}a_{\xi_5}\overline{e_{\xi_6}}}{\phi(\xi_1,\xi_2,\xi_3,\xi_4,\xi_5,\xi_6)}e^{-it\phi(\xi_1,\xi_2,\xi_3,\xi_4,\xi_5,\xi_6)},
\end{eqnarray}
where $\delta\gg 1$ is a fixed large constant, in fact, as like $\delta=10^{10}$.

First, we prepare the following lemma to estimate the convolution sum.
\begin{lemma}\label{lem:quintic}
Let $s>1$, and let $(c_{j,l})_{l\in 2\pi\mathbb{Z}}~(1\le j\le 6)$ be sequences of non-negative numbers.
Then it follows that
\begin{eqnarray}\label{qui-1}
\sum_{\xi\in 2\pi\mathbb{Z}/L}c_{j,\xi}\lesssim L^{1/2}\left(\sum_{\xi\in 2\pi \mathbb{Z}/L}m(\xi)^2c_{j,\xi}^2\right)^{1/2},
\end{eqnarray}
and
\begin{eqnarray}\label{qui-2}
\sum_{\xi_6\in 2\pi \mathbb{Z}/L}\sum^*_{m(\xi_6)\lesssim \max\{m(\xi_l)\mid 1\le l\le 5\}}
m(\xi_6)^2\prod_{j=1}^6c_{j,\xi_j}\lesssim L^2\prod_{j=1}^6\left(\sum_{\xi\in 2\pi \mathbb{Z}/L}m(\xi)^2c_{j,\xi}^2\right)^{1/2},
\end{eqnarray}
where constants term to the right-hand side are independent of $k$.
\end{lemma}
\noindent
{\it Proof.}
The proof of (\ref{qui-2}) can be obtained by considering (\ref{qui-1}).
The inequality (\ref{qui-1}) follows from the following fact
$$
\sum_{\xi\in 2\pi\mathbb{Z}}\frac{1}{m(\xi)^2}\lesssim 1.
$$
This is because a straightforward calculation as
\begin{eqnarray*}
\sum_{\eta=-100}^{100}\sum_{|\tau|\le k/2}\frac{1}{\langle \tau\rangle^{2(s-1/2)}}+\sum_{|\eta|\ge 100}\sum_{\tau=0}^{k-1}\frac{1}{\langle k\eta+\tau\rangle^{2s}}
 \lesssim  1+\int_{\mathbb{R}}\frac{dt}{\langle t\rangle^{2s}}\lesssim 1
\end{eqnarray*}
for every $s>1$.
\qed

\begin{lemma}\label{near-R}
Let ${\cal A}_r$ be
$$
{\cal A}_r=\left\{\xi=2\pi\left(k\widetilde{\eta}+\widetilde{\tau}+\frac{j}{L}\right)\mid \widetilde{\eta}\in\{0,1,3,4\},~(\widetilde{\tau},j)\in\mathbb{Z}^2,~|\widetilde{\tau}|\lesssim 1,~j\in[0,L)\right\}.
$$
Assume that $e(t)=(e_{\xi}(t))_{\xi\in 2\pi\mathbb{Z}/L}$ satisfies $e_{\xi}(0)=0$ for $\xi\in 2\pi\mathbb{Z}/L$ and
$$
\sup_{0\le t\le T}\|e(t)\|_s\lesssim \nu^{1/2}.
$$
Then for $|t|\le T$,
\begin{eqnarray}\label{near-R1}
\left(\sum_{\xi\in {\cal A}_r}m(\xi)^2|e_{\xi}(t)|^2\right)^{1/2}\lesssim \nu^{5/2}t.
\end{eqnarray}
and
\begin{eqnarray}\label{near-R2}
\int_0^t\left(\sum_{\xi\in {\cal A}_r}m(\xi)^2|e_{\xi}(t')|^2\right)^{1/2}\,dt'\lesssim \nu^{5/2}t^2.
\end{eqnarray}
\end{lemma}
\noindent
{\it Proof.}
Denote
$$
F(t)=\left(\int_0^t\sum_{\xi\in {\cal A}_r}m(\xi)^2|e_{\xi}(t')|^2\,dt'\right)^{1/2}.
$$
We note that $m(\xi)\sim 1$ for $\xi\in {\cal A}_r$.
From the assumption, equations (\ref{rfnls}), (\ref{RFNLS}) and Lemma \ref{lem:quintic}, we see that
\begin{eqnarray}\label{near-R3}
\sum_{\xi\in {\cal A}_r}m(\xi)^2|e_{\xi}(t)|^2\lesssim \nu^{5/2}\int_0^t\left(\sum_{\xi\in {\cal A}_r}m(\xi)^2|e_{\xi}(t')|^2\right)^{1/2}\,dt'
\le \nu^{5/2}t^{1/2}F(T_1)
\end{eqnarray}
for $t\le T_1\le T$.
Performing the integration on the time interval $[0,T_1]$, we have
$$
F(T_1)^2\lesssim \nu^{5/2}T_1^{3/2}F(T_1),
$$
which yields $F(T_1)\lesssim \nu^{5/2}T_1^{3/2}$ so that $F(t)\lesssim \nu^{5/2}t^{3/2}$ for $0\le t\le T$.
Inserting into (\ref{near-R3}) gives (\ref{near-R1}) and
$$
\int_0^t\left(\sum_{\xi\in {\cal A}_r}m(\xi)^2|e_{\xi}(t')|^2\right)^{1/2}\,dt'\lesssim t^{1/2}\nu^{5/2}t^{3/2}=\nu^{5/2}t^2,
$$
which implies (\ref{near-R2}).
\qed

\begin{lemma}\label{lem:men}
Let $\nu>0$ be a small constant.
Then
$$
\left|\widetilde{E}(t)-\|e(t)\|_s^2\right|\ll  \|e(t)\|_s^2.
$$
\end{lemma}
\noindent
{\it Proof.} 
It suffices to show that
\begin{eqnarray}\label{men-1}
\left|\frac{1}{L^5}\sum_{\xi_6}\sum_{\scriptstyle (\xi_1,\xi_2,\xi_3,\xi_4,\xi_5,\xi_6)\in A_1 \atop{\scriptstyle |\phi(\xi_1,\xi_2,\xi_3,\xi_4,\xi_5,\xi_6)|>\delta}}^*m(\xi_6)^2\frac{a_{\xi_1}\overline{a_{\xi_2}}a_{\xi_3}\overline{a_{\xi_4}}a_{\xi_5}\overline{e_{\xi_6}}}{\phi(\xi_1,\xi_2,\xi_3,\xi_4,\xi_5,\xi_6)}e^{-it\phi(\xi_1,\xi_2,\xi_3,\xi_4,\xi_5,\xi_6)}\right| \ll  \|e(t)\|_s^2.
\end{eqnarray}

In the case when $\xi_6\in {\cal A}_r$, we distinguish two cases:
\begin{itemize}
\item
$\xi_m\in{\cal R}$ for all $1\le m\le 5$ in the sum (\ref{men-1}),
\item
there is at least one element of $\xi_m\in{\cal R}~(1\le m\le 5)$ such that $\xi_m\not\in{\cal R}$ in the sum (\ref{men-1}). 
\end{itemize} 
Consider the first case.
We use the formula in (\ref{tilde}) such that for $1\le m\le 6$,
$$
\xi_m=2\pi\left(k\widetilde{\eta}_m+\widetilde{\tau}_m+\frac{j_m}{L}\right),
$$
where $\widetilde{\eta}_1=\widetilde{\eta}_2=\widetilde{\eta}_3=\widetilde{\eta}_4=\widetilde{\eta}_5=\widetilde{\eta}_6\in\{0,1,3,4\},~\widetilde{\tau}_m=0~(1\le m\le 5)$ and $|\widetilde{\tau}_6|\lesssim 1$, since by $(\xi_1,\xi_2,\xi_3,\xi_4,\xi_5,\xi_6)\in A_1$ and $\xi_6\in{\cal A}_r$.
Then we obtain
\begin{eqnarray}\label{R1-11}
\left(\widetilde{\tau}_1+\frac{j_1}{L}\right)-\left(\widetilde{\tau}_2+\frac{j_2}{L}\right)+\left(\widetilde{\tau}_3+\frac{j_3}{L}\right)-\left(\widetilde{\tau}_4+\frac{j_4}{L}\right)+\left(\widetilde{\tau}_5+\frac{j_5}{L}\right)-\left(\widetilde{\tau}_6+\frac{j_6}{L}\right)=0
\end{eqnarray}
and
\begin{eqnarray}
& & \phi(\xi_1,\xi_2,\xi_3,\xi_4,\xi_5,\xi_6) \nonumber\\
& = &\left(\widetilde{\tau}_1+\frac{j_1}{L}\right)^2-\left(\widetilde{\tau}_2+\frac{j_2}{L}\right)^2+\left(\widetilde{\tau}_3+\frac{j_3}{L}\right)^2-\left(\widetilde{\tau}_4+\frac{j_4}{L}\right)^2+\left(\widetilde{\tau}_5+\frac{j_5}{L}\right)^2-\left(\widetilde{\tau}_6+\frac{j_6}{L}\right)^2.\label{R1-12}
\end{eqnarray}
In the case when $\xi_m\in{\cal R}$ for all $1\le m\le 5$ in the sum (\ref{men-1}), we have $|\phi(\xi_1,\xi_2,\xi_3,\xi_4,\xi_5,\xi_6)|\lesssim 1$, which is out of range in the sum.

In the second case, from $a_{\xi_m}=e_{\xi_m}$ for $\xi\not\in{\cal R}$, the contribution of this case to the left-hand side of (\ref{men-1}) is bounded by
\begin{eqnarray}\label{R1-13}
c\frac{\nu^2}{\delta}\|e(t)\|_s^2\ll  \|e(t)\|_s^2.
\end{eqnarray}

Next we may suppose $\xi_6\not\in {\cal A}_r$ so that there is at least one element of $\xi_m~(1\le m\le 5)$ such that $\xi_m\not\in {\cal R},~a_{\xi_m}=e_{\xi_m}$ and $|\widetilde{\tau}_6|\lesssim |\widetilde{\tau}_m|$.
In this case, it is easy see that $m(\xi_6)\lesssim m(\xi_m)$ so that the contribution of this case to the left-hand side of (\ref{men-1}) is bounded by the same bound as in (\ref{R1-13}), 
which is acceptable.
%
\qed

A crucial step in the proof of Theorem \ref{thm-main} is to establish the following proposition.
\begin{proposition}\label{prop:approx}
Let $s\in(1,3/2]$ and $T>0$.
Given $a(0)=(a_{\xi}(0))_{\xi\in 2\pi\mathbb{Z}/L}\in \ell^2$ and $r(0)=(r_{\xi}(0))_{\xi\in {\cal R}}$, let $e(t)=(e_{\xi}(t))_{\xi\in 2\pi\mathbb{Z}/L}$ be a solution to the equation (\ref{eq:e}) with initial data $e(0)=(e_{\xi}(0))_{\xi\in 2\pi\mathbb{Z}/L}$ for $e_{\xi}(0)=a_{\xi}(0)-r_{\xi}(0)=0$ on times $0\le t\le T$.
Assume that
$$
\sup_{0\le t\le T}\|e(t)\|_s^2\lesssim \nu.
$$
Then
\begin{eqnarray}\label{app2}
\|e(t)\|_s^2\lesssim \frac{\nu^3}{\langle k\rangle^{5/2-s}}e^{c\nu^2 t}+\frac{\nu^3}{\langle k\rangle^{3/2-s}}\left(e^{c\nu^2 t}-1\right)+\nu\left(e^{c\nu^2t}-1-c\nu^2t\right).
\end{eqnarray}
\end{proposition}
\noindent
{\it Proof of Proposition \ref{prop:approx}.}
We notice $\widetilde{E}(0)=0$.
By the fundamental theorem of calculus, we see that from (\ref{eq:e})
\begin{eqnarray}
\widetilde{E}(t)= -\frac{2}{L}\Im\sum_{n=1}^3\sum_{\xi\in 2\pi\mathbb{Z}/L}\int_0^t m(\xi)^2\widetilde{R}_{\xi}^n(t')\overline{e_{\xi}(t')} \,dt'+\widetilde{R}_4(t)\label{dif}
\end{eqnarray}
where
\begin{eqnarray*}
\widetilde{R}_{\xi}^1 (t)= \frac{1}{L^4}\sum_{\scriptstyle (\xi_1,\xi_2,\xi_3,\xi_4,\xi_5,\xi)\in A_1 \atop{\scriptstyle |\phi(\xi_1,\xi_2,\xi_3,\xi_4,\xi_5,\xi)|\le \delta}}^*
a_{\xi_1}(t)\overline{a_{\xi_2}(t)}a_{\xi_3}(t)\overline{a_{\xi_4}(t)}a_{\xi_5}(t)e^{-it\phi(\xi_1,\xi_2,\xi_3,\xi_4,\xi_5,\xi)},
\end{eqnarray*}
\begin{eqnarray*}
\widetilde{R}_4(t)=\frac{2}{L^5}\Re \sum_{\xi_6}\sum_{\scriptstyle (\xi_1,\xi_2,\xi_3,\xi_4,\xi_5,\xi_6)\in A_1 \atop{\scriptstyle |\phi(\xi_1,\xi_2,\xi_3,\xi_4,\xi_5,\xi_6)|>\delta}}^*\int_0^t m(\xi_6)^2\frac{\partial_{t'}\left(a_{\xi_1}(t')\overline{a_{\xi_2}(t')}a_{\xi_3}(t')\overline{a_{\xi_4}(t')}a_{\xi_5}(t')\overline{e_{\xi_6}(t')}\right)}{\phi(\xi_1,\xi_2,\xi_3,\xi_4,\xi_5,\xi_6)}e^{-it'\phi(\xi_1,\xi_2,\xi_3,\xi_4,\xi_5,\xi_6)}\,dt'
\end{eqnarray*}
and $\widetilde{R}_{\xi}^n(t)=R_{\xi}^n(t)$ for $n=2,3$.

For the terms $\widetilde{R}_{\xi}^n~(n=1,2)$, we have the following two lemmas.
\begin{lemma}\label{lem:R1}
\begin{eqnarray*}
\frac{1}{L}\sum_{\xi\in 2\pi\mathbb{Z}/L}\int_0^t m(\xi)^2\left|\widetilde{R}_{\xi}^1(t')\overline{e_{\xi}(t')}\right| \,dt'\lesssim \left(\delta^{1/4}+\frac{1}{L^{1/2}}\right)\nu^5t^2+\nu^2\int_0^t \|e(t')\|_s^2\,dt'.
\end{eqnarray*}
\end{lemma}
\noindent
{\it Proof.}
It suffices to show that
\begin{eqnarray}
& & \frac{1}{L^5}\sum_{\xi_6}\sum_{\scriptstyle (\xi_1,\xi_2,\xi_3,\xi_4,\xi_5,\xi_6)\in A_1 \atop{\scriptstyle |\phi(\xi_1,\xi_2,\xi_3,\xi_4,\xi_5,\xi_6)|\le \delta}}^*
m(\xi_6)^2|a_{\xi_1}(t)\overline{a_{\xi_2}(t)}a_{\xi_3}(t)\overline{a_{\xi_4}(t)}a_{\xi_5}(t)\overline{e_{\xi}(t')}|\,dt'\nonumber \\
& \lesssim& \left(\delta^{1/2}+\frac{1}{L}\right)^{1/2}\nu^5t^2+\nu^2\int_0^t \|e(t')\|_s^2\,dt'.\label{R1-1}
\end{eqnarray}

We note that if $(\xi_1,\xi_2,\xi_3,\xi_4,\xi_5,\xi_6)\in A_1$, then any of $\widetilde{\eta}_m~(1\le m\le 6)$ coincide with each other, and the form of $m(\xi_m)$ is written $m(\xi_m)=\langle\widetilde{\tau}_m\rangle^{s-1/2}$.
By using the notation in (\ref{tilde}), we have (\ref{R1-11}) and (\ref{R1-12}).
The restriction $|\phi(\xi_1,\xi_2,\xi_3,\xi_4,\xi_5,\xi_6)|< \delta$ implies that for fixed four elements of $\xi_m~(1\le m\le 6)$ under the plane $\xi_1+\xi_3+\xi_5=\xi_2+\xi_4+\xi_6$, the cardinalities of one of the other elements is at most $c\delta^{1/2}L+1$.

In the case when $m(\xi_6)\lesssim 1$ in the sum of (\ref{R1-1}), by the above observation, we have the bound of the left-hand side of (\ref{R1-1}) by
$$
c\left(\frac{\delta^{1/2}L+1}{L}\right)^{1/2}\nu^{5/2}\int_0^t\|e(t')\|_s\,dt',
$$
which is bounded by
$$
c\left(\delta^{1/2}+\frac{1}{L}\right)^{1/2}\nu^5t^2,
$$
since by Lemma \ref{near-R}.

In the case when $m(\xi_6)\gg 1$, from (\ref{R1-11}), we have
$$
1\ll m(\xi_6)\lesssim \max_{1\le m\le 5}m(\xi_m),
$$
which implies that at least there exists one element $\xi_0$ of $\xi_m~(1\le m\le 5)$ in the sum such that $\xi_0\not\in{\cal R}$ and $m(\xi_6)\sim m(\xi_0)$.
Then the contribution of this case to the left-hand side of (\ref{R1-1}) is bounded by
$$
c\nu^2\int_0^t\|e(t')\|_s^2\,dt'.
$$
Thus the proof of that result is complete.
\qed

\begin{lemma}\label{lem:R2}
\begin{eqnarray*}
\frac{1}{L}\sum_{\xi\in 2\pi\mathbb{Z}/L}\int_0^t m(\xi)^2\left|\widetilde{R}_{\xi}^2(t')\overline{e_{\xi}(t')}\right| \,dt'\lesssim \nu^2\int_0^t \|e(t')\|_s^2\,dt'.
\end{eqnarray*}
\end{lemma}
\noindent
{\it Proof.}
The proof of this lemma follow easily from the form of $R_{\xi}^2$ and the fact that $R_{\xi}^2=0$ for $\xi\not\in{\cal R}$, so that shall be omitted.
\qed

To bound the last factor on the right-hand side of (\ref{dif}), we have the following.
\begin{lemma}\label{lem:R4}
\begin{eqnarray*}
|\widetilde{R}_4(t)|\lesssim \frac{\nu^4}{\delta}\int_0^t\|e(t')\|_s^2\,dt'.
\end{eqnarray*}
\end{lemma}
\noindent
{\it Proof.}
The proof of this lemma is a straightforward.
The main objective is to use equations (\ref{rfnls}), (\ref{RFNLS}) and the a priori bound.

The contribution to the upper bound in term of $|\widetilde{R}_4(t)|$ is the sums of the following two terms:
\begin{eqnarray*}
\frac{1}{\delta L^5}\sum_{\xi_6}\sum_{\scriptstyle (\xi_1,\xi_2,\xi_3,\xi_4,\xi_5,\xi_6)\in A_1 \atop{\scriptstyle |\phi(\xi_1,\xi_2,\xi_3,\xi_4,\xi_5,\xi_6)|>\delta}}^*\int_0^t m(\xi_6)^2\left|\partial_{t'}\left(a_{\xi_1}(t')\overline{a_{\xi_2}(t')}a_{\xi_3}(t')\overline{a_{\xi_4}(t')}a_{\xi_5}(t')\right)\right|\left|\overline{e_{\xi_6}(t))}\right|\,dt'
\end{eqnarray*}
and
\begin{eqnarray*}
\frac{1}{\delta L^5}\sum_{\xi_6}\sum_{\scriptstyle (\xi_1,\xi_2,\xi_3,\xi_4,\xi_5,\xi_6)\in A_1 \atop{\scriptstyle |\phi(\xi_1,\xi_2,\xi_3,\xi_4,\xi_5,\xi_6)|>\delta}}^*\int_0^t m(\xi_6)^2\left|a_{\xi_1}(t')\overline{a_{\xi_2}(t')}a_{\xi_3}(t')\overline{a_{\xi_4}(t')}a_{\xi_5}(t')\right|\left|\partial_{t'}\overline{e_{\xi_6}(t)}\right|\,dt'
\end{eqnarray*}
Here we only give a proof for the first term, because the second term can be handled similarly. 

In the same way as in the proof of Lemma \ref{lem:men}, we use again the formula in (\ref{tilde}), where $\widetilde{\eta}_1=\widetilde{\eta}_2=\widetilde{\eta}_3=\widetilde{\eta}_4=\widetilde{\eta}_5=\widetilde{\eta}_6\in\{0,1,3,4\}$ and $-k/2<\widetilde{\tau}_m\le k/2$ for $1\le m\le 6$.

Notice that by (\ref{R1-11}) and (\ref{R1-12}), the restriction $|\phi(\xi_1,\xi_2,\xi_3,\xi_4,\xi_5,\xi_6)|>\delta$ in the sum $\widetilde{R}_4(t)$ yields that at least two elements of $\xi_m~(1\le m\le 6)$ satisfy $|\widetilde{\tau}_m|\gg 1$, in which we denote $\xi_a,~\xi_b$, namely, $m(\xi_a)\gg 1,~m(\xi_b)\gg 1$.
We may assume $m(\xi_6)\lesssim \min\{m(\xi_a),m(\xi_b)\}$.

Observe that if $|\widetilde{\tau}_a|\gg 1$, then by (\ref{rfnls}), $e_{\xi_a}(t)=a_{\xi_a}(t)$ satisfies the equation:
\begin{eqnarray}\label{eq:xi_a}
i\dot{a}_{\xi_a}=\frac{1}{L^4}\sum^*a_{\xi_1}\overline{a_{\xi_2}}a_{\xi_3}\overline{a_{\xi_4}}a_{\xi_5}e^{-it\phi(\xi_1,\xi_2,\xi_3,\xi_4,\xi_5,\xi_a)}.
\end{eqnarray}
There is at least one $a_{\xi_m}~(1\le m\le 5)$ of each quintic nonlinearity $a_{\xi_1}\overline{a_{\xi_2}}a_{\xi_3}\overline{a_{\xi_4}}a_{\xi_5}$ in (\ref{eq:xi_a}) such that $m(\xi_m)\gtrsim m(\xi_a)$.
This is easy to prove by writing $\xi_m$ as the formula (\ref{tilde}).

In view of above, in the case that the derivative $\partial_{t'}$ falls into either the function associated to $\xi_a$ or $\xi_b$, we have the bound of $|\widetilde{R}_4(t)|$ by
$$
c\frac{\nu^4}{\delta}\int_0^t\|e(t')\|_s^2\,dt'.
$$
On the other hand, in the case that the derivative $\partial_{t'}$ falls into the function associated with neither $\xi_a$ or $\xi_b$, we have the same bound as above.

Thus the proof is complete.
\qed

It remains to estimate the term
\begin{eqnarray}\label{R3}
-\frac{2}{L}\Im\sum_{\xi\in 2\pi\mathbb{Z}/L}\int_0^t m(\xi)^2\widetilde{R}_{\xi}^3(t')\overline{e_{\xi}(t')} \,dt'.
\end{eqnarray}

\subsection{Several technical lemmas}
\noindent

First we shall prepare several lemmas.
We keep the convention of notation expressed at Definition \ref{def:tilde}.
For frequencies $\xi\in 2\pi\mathbb{Z}/L$, we use the expression in (\ref{expr:xi}) and (\ref{tilde}).

\begin{lemma}\label{rc-1}
Let $\xi_n\in {\cal R}~(1\le n\le 5)$ and $\xi_6\in {\cal R}^c$ be elements of expression as described in (\ref{expr:xi}).
If $|\phi(\xi_1,\xi_2,\xi_3,\xi_4,\xi_5,\xi_6)|\ll\langle k\rangle^2$, then $\widetilde{\eta}_6\in\{0,1,3,4\}$ and $\widetilde{\tau}_6\in\{\pm 1,\pm 2\}$.
\end{lemma}
\noindent
{\it Proof.}
The proof is straightforward.
Expressing $\xi_n$ as described in (\ref{expr:xi}) and inserting into $\xi_1+\xi_3+\xi_5=\xi_2+\xi_4+\xi_6$, we see that
$$
\eta_n\in\{0,\,1,\,3,\,4\},~\tau_n=0\quad\mbox{for $1\le n\le 5$},
$$
and
$$
\eta_1+\eta_3+\eta_5-\eta_2-\eta_4-\eta_6=\frac{j_2+j_4+j_6-j_1-j_3-j_5}{Lk}+\frac{\tau_6}{k}.
$$
We analyze and discuss each of the case study with respect to $\tau_6$.

In the case when $3\le\tau_6\le k-3$, clearly
$$
\eta_1+\eta_3+\eta_5-\eta_2-\eta_4-\eta_6\in (0,1),
$$
which is a contradiction, since the left-hand side is integer number.
So the remaining cases are $0\le \tau_6\le 2$, and $k-2\le \tau_6\le k-1$.

Consider the case when $0\le \tau_6\le  2$.
Clearly
\begin{eqnarray*}
\eta_1+\eta_3+\eta_5=\eta_2+\eta_4+\eta_6+O\left(\frac{1}{k}\right),
\end{eqnarray*}
\begin{eqnarray*}
(0,1)\ni \frac{1}{4\pi^2 k^2}\phi(\xi_1,\xi_2,\xi_3,\xi_4,\xi_5,\xi_6)=\eta_1^2-\eta_2^2+\eta_3^2-\eta_4^2+\eta_5^2-\eta_6^2+O\left(\frac{1}{k}\right)
\end{eqnarray*}
and hence
\begin{eqnarray*}
\eta_1+\eta_3+\eta_5 =\eta_2+\eta_4+\eta_6,\quad
\eta_1^2+\eta_3^2+\eta_5^2 = \eta_2^2+\eta_4^2+\eta_6^2.
\end{eqnarray*}
We distinguish the cases.

In the case when one element of $\eta_1,\,\eta_3,\,\eta_5$ is equal to one of $\eta_2,\,\eta_4$, we may easily obtain $\{\eta_1,\,\eta_3,\,\eta_5\}=\{\eta_2,\,\eta_4,\,\eta_6\}$, which leads to $\eta_6\in\{0,\,1,\,3,\,4\}$.

Next we consider the case when any elements of $\eta_1,\,\eta_3,\,\eta_5$ do not match any of $\eta_2,\,\eta_4$, namely $\eta_1,\,\eta_3,\,\eta_5\not\in\{\eta_2,\eta_4\}$.
When two elements of $\eta_1,\,\eta_3,\,\eta_5$ are to be considered equal, we may assume $\eta_1=\eta_3$ by a symmetry argument.
Observe that $\eta_6$ satisfies the one of
\begin{eqnarray}\label{eta_6}
\frac{1}{2}\left(2\eta_1+\eta_5-\eta_2\pm\sqrt{(\eta_2-\eta_5)(4\eta_1-\eta_5-3\eta_2)}\right),
\end{eqnarray}
which should take in an integer number.
We calculate the considering numerical value of (\ref{eta_6}) in all cases.
In the case when $\eta_1=\eta_5$, we have that $\eta_1$ and $\eta_2$ satisfy $\eta_1=\eta_2=\eta_4$, which is not acceptable.
In the case when $\eta_1=\eta_3\ne\eta_5$, we will check all cases under the condition $\eta_1,\,\eta_3,\,\eta_5\not\in\{\eta_2,\eta_4\}$: 
\begin{itemize}
\item
$\{\eta_1,\,\eta_5\}=\{0,\,1\},~\mbox{then}~\sqrt{(\eta_2-\eta_5)(4\eta_1-\eta_5-3\eta_2)}\not\in\mathbb{Z}~\mbox{for}~\eta_2\in\{3,\,4\}$,
\item
$\{\eta_1,\,\eta_5\}=\{0,\,3\},~\mbox{then}~\sqrt{(\eta_2-\eta_5)(4\eta_1-\eta_5-3\eta_2)}\in\mathbb{Z}~\mbox{if and only if}~\eta_1=3,\,\eta_5=0,\,\eta_2\in\{1,\,4\}$,
\item
$\{\eta_1,\,\eta_5\}=\{0,\,4\},~\mbox{then}~\sqrt{(\eta_2-\eta_5)(4\eta_1-\eta_5-3\eta_2)}\in\mathbb{Z}~\mbox{for}~\eta_2\in\{1,\,3\}$,
\item
$\{\eta_1,\,\eta_5\}=\{1,\,3\},~\mbox{then}~\sqrt{(\eta_2-\eta_5)(4\eta_1-\eta_5-3\eta_2)}\not\in\mathbb{Z}~\mbox{for}~\eta_2\in\{0,\,4\}$,
\item
$\{\eta_1,\,\eta_5\}=\{1,\,4\},~\mbox{then}~\sqrt{(\eta_2-\eta_5)(4\eta_1-\eta_5-3\eta_2)}\in\mathbb{Z}~\mbox{if and only if}~\eta_1=1,\,\eta_5=4,\,\eta_2\in\{0,\,3\}$,
\item
$\{\eta_1,\,\eta_5\}=\{3,\,4\},~\mbox{then}~\sqrt{(\eta_2-\eta_5)(4\eta_1-\eta_5-3\eta_2)}\not\in\mathbb{Z}~\mbox{for}~\eta_2\in\{0,\,1\}$.
\end{itemize}
Therefore there are two possible cases:
\begin{itemize}
\item
two elements of $\eta_1,\,\eta_3,\,\eta_5$ are $1$ and the rest of those elements is $4$, furthermore two elements of $\eta_2,\,\eta_4,\,\eta_6$ are $3$ and the rest of those elements is $0$, 
\item
two elements of $\eta_2,\,\eta_4,\,\eta_6$ are $1$ and the rest of those elements is $4$, moreover two elements of $\eta_1,\,\eta_3,\,\eta_5$ are $3$ and the rest of those elements is $0$.
\end{itemize}
We remark that $\tau_6\ne 0$ since by $\xi_6\not\in{\cal R}$, that are acceptable.

In the case when any two elements of $\eta_1,\,\eta_3,\,\eta_5$ do not match with each other, we again consider all cases:
\begin{itemize}
\item
$\{\eta_1,\,\eta_3,\,\eta_5\}=\{0,\,1,\,3\},~\eta_2=\eta_4=4,~\eta_6=-4$,
\item
$\{\eta_1,\,\eta_3,\,\eta_5\}=\{1,\,3,\,4\},~\eta_2=\eta_4=0,~\eta_6=8$,
\item
$\{\eta_1,\,\eta_3,\,\eta_5\}=\{0,\,3,\,4\},~\eta_2=\eta_4=1,~\eta_6=5$,
\item
$\{\eta_1,\,\eta_3,\,\eta_5\}=\{0,\,1,\,4\},~\eta_2=\eta_4=3,~\eta_6=-1$,
\end{itemize}
which are inadequate for the equation $\eta_1^2+\eta_3^2+\eta_5^2 = \eta_2^2+\eta_4^2+\eta_6^2$.

Consider next the case when $k-2\le \tau_6\le k-1$.
Obviously
\begin{eqnarray*}
\eta_1+\eta_3+\eta_5=\eta_2+\eta_4+\widetilde{\eta}_6+O\left(\frac{1}{k}\right),
\end{eqnarray*}
\begin{eqnarray*}
(0,1)\ni\frac{1}{4\pi^2 k^2}\phi(\xi_1,\xi_2,\xi_3,\xi_4,\xi_5,\xi_6)=\eta_1^2-\eta_2^2+\eta_3^2-\eta_4^2+\eta_5^2-\widetilde{\eta}_6^2+O\left(\frac{1}{k}\right)
\end{eqnarray*}
and hence
\begin{eqnarray*}
\eta_1+\eta_3+\eta_5 = \eta_2+\eta_4+\widetilde{\eta}_6,\quad
\eta_1^2+\eta_3^2+\eta_5^2  =  \eta_2^2+\eta_4^2+\widetilde{\eta}_6^2.
\end{eqnarray*}
Using the terminology in the previous cases, we have that $\widetilde{\eta}_6\in\{0,\,1,\,3,\,4\}$.

This concludes the proof of the lemma.
\qed
%
%

\begin{corollary}\label{s-1}
Let $\eta_j\in\{0,\,1,\,3,\,4\}~(j=1,\,2,\,3,\,5)$, and let $\eta_4,\,\eta_6\in\mathbb{Z}$ satisfy
$$
\eta_1+\eta_3+\eta_5=\eta_2+\eta_4+\eta_6,\quad
\eta_1^2+\eta_3^2+\eta_5^2=\eta_2^2+\eta_4^2+\eta_6^2.
$$
Then either one of the following holds:
\begin{itemize}
\item
$\{\eta_1,\,\eta_3,\,\eta_5\}=\{\eta_2,\,\eta_4,\,\eta_6\}$,
\item
two elements of $\{\eta_1,\,\eta_3,\,\eta_5\}$ are $1$ and the rest of those elements is $4$; two elements of $\{\eta_2,\,\eta_4,\eta_6\}$ are $3$ and the rest of those elements is $0$,
\item
symmetric case of above replacing $\{\eta_1,\,\eta_3,\,\eta_5\}$ with $\{\eta_2,\,\eta_4,\eta_6\}$ and $\{\eta_2,\,\eta_4,\eta_6\}$ with  $\{\eta_1,\,\eta_3,\,\eta_5\}$.
\end{itemize}
\end{corollary}
\noindent
{\it Proof.}
In the case when any two elements of $\eta_1,\,\eta_3,\,\eta_5$ do not match with each other, we consider all cases:
\begin{itemize}
\item
$\{\eta_1,\,\eta_3,\,\eta_5\}=\{0,\,1,\,3\},~\eta_2=4,~\eta_4+\eta_6=0$,
\item
$\{\eta_1,\,\eta_3,\,\eta_5\}=\{1,\,3,\,4\},~\eta_2=0,~\eta_4+\eta_6=8$,
\item
$\{\eta_1,\,\eta_3,\,\eta_5\}=\{0,\,3,\,4\},~\eta_2=1,~\eta_4+\eta_6=6$,
\item
$\{\eta_1,\,\eta_3,\,\eta_5\}=\{0,\,1,\,4\},~\eta_2=3,~\eta_4+\eta_6=2$,
\end{itemize}
however $\eta_1^2+\eta_3^2+\eta_5^2 = \eta_2^2+\eta_4^2+\eta_6^2$ is no longer satisfied.
Therefore, the proof of Lemma \ref{rc-1} permits us to conclude the proof of the corollary.
\qed

\subsection{Contribution of $\xi\in {\cal A}_r$ to (\ref{R3})}\label{sec:A_r}
\noindent
 
Let us first consider the contribution of $\xi\in{\cal A}_r$ to (\ref{R3}).
Note that $m(\xi)\lesssim 1$ for $\xi\in{\cal A}_r$.
%
%
The contribution of this case to (\ref{R3}) is
\begin{eqnarray}\label{R-3}
-\frac{2}{L^5}\Im \sum_{\xi_6\in{\cal A}_r}\sum^*_{\scriptstyle \text{res}(\xi_6)^c \atop{\scriptstyle (\xi_1,\xi_2,\xi_3,\xi_4,\xi_5,\xi_6)\in A_1^c}}\int_0^t m(\xi_6)^2a_{\xi_1}(t')\overline{a_{\xi_2}(t')}a_{\xi_3}(t')\overline{a_{\xi_5}(t')} a_{\xi_5}(t') \overline{e_{\xi_6}(t')} e^{-it'\phi(\xi_1,\xi_2,\xi_3,\xi_4,\xi_5,\xi_6)}\,dt'.
\end{eqnarray}
Note that $m(\xi)\lesssim 1$ for $\xi\in{\cal A}_r$.
By Lemma \ref{near-R}, we have that the contribution of this case is bounded by
\begin{eqnarray}
c\nu^{5/2}\int_0^t \|e(t')\|_s\,dt'\lesssim \nu^5t^2.
\end{eqnarray}

\subsection{Contribution of $\xi\not\in {\cal A}_r$ to (\ref{R3})}
\noindent

%
%
Consider the associated function
\begin{eqnarray}\label{rem-4}
\frac{1}{L^5}\Im \int_0^t\sum_{\xi_6\not\in{\cal A}_r}\sum^*_{\scriptstyle \text{res}(\xi_6)^c \atop{\scriptstyle (\xi_1,\xi_2,\xi_3,\xi_4,\xi_5,\xi_6)\in A_1^c}}m(\xi_6)^2a_{\xi_1}(t')\overline{a_{\xi_2}(t')}a_{\xi_3}(t')\overline{a_{\xi_4}(t')}a_{\xi_5}(t')\overline{e_{\xi_6}(t')}e^{-it'\phi(\xi_1,\xi_2,\xi_3,\xi_4,\xi_5,\xi_6)}\,dt'.
\end{eqnarray}
We split the sums into three cases:
\begin{itemize}
\item[(i)]
$\max_{1\le m\le 5}m(\xi_m)\gtrsim m(\xi_6)$,
\item[(ii)]
$\max_{1\le m\le 5}m(\xi_m)\ll m(\xi_6)$; and at least four elements of $\xi_m~(1\le m\le 5)$ satisfy $m(\xi_m)\ll \langle k\rangle^{s-1/2}$,
\item[(iii)]
$\max_{1\le m\le 5}m(\xi_m)\ll m(\xi_6)$; and at least two elements of $\xi_m~(1\le m\le 5)$ satisfy $m(\xi_m)\gtrsim \langle k\rangle^{s-1/2}$.
\end{itemize}

\underline{Case (i)}:
In this case, we easily see that this contribution to (\ref{rem-4}) is bounded by
\begin{eqnarray}\label{sq}
c\nu^2\int_0^t\|e(t')\|_s^2\,dt'.
\end{eqnarray}

\underline{Case (ii)}:
In this case, we formulate $\xi_j~(1\le j\le 6)$ as in the form (\ref{tilde}) such that
\begin{eqnarray}\label{441}
\xi_m=2\pi\left(k\widetilde{\eta}_m+\widetilde{\tau}_m+\frac{j_m}{L}\right).
\end{eqnarray}
and distinguish the following two cases:
\begin{itemize}
\item[(ii-1)]
$|\widetilde{\eta}_6|>50$,
\item[(ii-2)]
$|\widetilde{\eta}_6|\le 50$.
\end{itemize}

In the case of (ii-1), we only need to work with $\widetilde{\eta}_m\in \{0,1,3,4\},~|\widetilde{\tau}_m|\ll k$ for $1\le m\le 4$.
The restriction that $\xi_1+\xi_3+\xi_5=\xi_2+\xi_4+\xi_6$ and $m(\xi_5)\ll m(\xi_6)$ yields $\widetilde{\eta}_5\le 100,~99\le \widetilde{\eta}_6,~0\le \widetilde{\eta}_6-\widetilde{\eta}_5<10,~ k^{s-1/2}\lesssim m(\xi_5)\ll m(\xi_6)\lesssim k^s$.

We first observe that in the case when $|\phi(\xi_1,\xi_2,\xi_3,\xi_4,\xi_5,\xi_6)|\gtrsim k^2$.
While integration by parts, time derivative hits either $a_{\xi_j}$ or $e_{\xi_6}$ in (\ref{rem-4}).
The contribution of this case to (\ref{rem-4}) is the sum of the following three terms
\begin{eqnarray*}
\frac{\langle k\rangle^{-3/2}}{L^5}\left|\sum_{\xi_6\not\in{\cal A}_r}\sum^*a_{\xi_1}(t)\overline{a_{\xi_2}(t)}a_{\xi_3}(t)\overline{a_{\xi_4}(t)}m(\xi_5) a_{\xi_5}(t) m(\xi_6)\overline{e_{\xi_6}(t)}\right|,
\end{eqnarray*}
\begin{eqnarray*}
\frac{\langle k\rangle^{-3/2}}{L^5}\int_0^t\sum_{\xi_6\not\in{\cal A}_r}\sum^*\left|\partial_{t'}\left(a_{\xi_1}(t')\overline{a_{\xi_2}(t')}a_{\xi_3}(t')\overline{a_{\xi_4}(t')} \langle k\rangle^{s-1/2}a_{\xi_5}(t')\right)m(\xi_6)\overline{e_{\xi_6}(t')}\right|\,dt'
\end{eqnarray*}
and
\begin{eqnarray*}
\frac{\langle k\rangle^{-1}}{L^5}\int_0^t\sum_{\xi_6\not\in{\cal A}_r}\sum^*\left|a_{\xi_1}(t')\overline{a_{\xi_2}(t')}a_{\xi_3}(t')\overline{a_{\xi_4}(t')} m(\xi_5)a_{\xi_5}(t')\langle k\rangle^{s-1/2}\partial_{t'}\overline{e_{\xi_6}(t')}\,dt'\right|,
\end{eqnarray*}
which are bounded by
\begin{eqnarray}
\lesssim \frac{\nu^2}{\langle k\rangle^{3/2}}\|e(t)\|_s^2+\nu^2\int_0^t\|e(t')\|_s^2\,dt',\label{rem4411}
\end{eqnarray}
for $s>1$.

Next suppose that $|\phi(\xi_1,\xi_2,\xi_3,\xi_4,\xi_5,\xi_6)|\ll k^2$ holds.
Observe
\begin{eqnarray}\label{lin}
k(\widetilde{\eta}_1-\widetilde{\eta}_2+\widetilde{\eta}_3-\widetilde{\eta}_4+\widetilde{\eta}_5-\widetilde{\eta}_6)-(\widetilde{\tau}_6-\widetilde{\tau}_5)=\frac{j_2+j_4+j_6-j_1-j_3-j_5}{L}\in(-3,3),
\end{eqnarray}
and
\begin{eqnarray}\label{442}
1\gg \frac{1}{(2\pi)^2}|\phi(\xi_1,\xi_2,\xi_3,\xi_4,\xi_5,\xi_6)|\ge \left|k(\widetilde{\eta}_5-\widetilde{\eta}_6)+\widetilde{\tau}_5-\widetilde{\tau}_6+\frac{j_5-j_6}{L}\right|(\xi_5+\xi_6)-50k^2.
\end{eqnarray}
Considering the case $m(\xi_5)\ll m(\xi_6)$, we investigate a possible case above and have two possible cases; that (a): $\widetilde{\eta}_6=\widetilde{\eta}_5=99$, and that (b): $\widetilde{\eta}_6=\widetilde{\eta}_5+1\in\{99,100\}$; moreover, subdivide each case into several cases as follows:
\begin{itemize}
\item[(a)]
$\widetilde{\eta}_6=\widetilde{\eta}_5=99,~\widetilde{\tau}_5<0,~\widetilde{\tau}_6\gg 1,~m(\xi_5)=\langle k\rangle^{s-1/2},~m(\xi_6)\sim \langle k\rangle^{s-1/2}\langle\widetilde{\tau}_6\rangle^{1/2}$,
\item[(b1)]
$\widetilde{\eta}_6=\widetilde{\eta}_5+1=99,~\widetilde{\tau}_6\gg 1,~m(\xi_5)=\langle k\rangle^{s-1/2},~m(\xi_6)=\langle k\rangle^{s-1/2}\langle\widetilde{\tau}_6\rangle^{1/2}$,
\item[(b2)]
$\widetilde{\eta}_6=\widetilde{\eta}_5+1=100,~\widetilde{\tau}_5<0,~m(\xi_5)=\langle k\rangle^{s-1/2},~m(\xi_6)\sim \langle k\rangle^{s}$,
\item[(b3)]
$\widetilde{\eta}_6=\widetilde{\eta}_5+1=100,~\widetilde{\tau}_5\ge 0,~\langle\widetilde{\tau}_5\rangle\ll\langle k\rangle,~m(\xi_5)=\langle k\rangle^{s-1/2}\langle\widetilde{\tau}_5\rangle^{1/2},~m(\xi_6)\sim \langle k\rangle^{s}$,
\end{itemize}

In the case of (a), we see that if $\widetilde{\tau}_6-\widetilde{\tau}_5<k-3$ then (\ref{lin}) implies
$$
1\ll \widetilde{\tau}_6<\widetilde{\tau}_6-\widetilde{\tau}_5<3,
$$
which is not subject to the considering case.
Then $|k(\widetilde{\eta}_5-\widetilde{\eta}_6)+\widetilde{\tau}_5-\widetilde{\tau}_6|\ge k-3$, however by (\ref{442}), it follows that
\begin{eqnarray}\label{count}
|\phi(\xi_1,\xi_2,\xi_3,\xi_4,\xi_5,\xi_6)|\gtrsim k^2,
\end{eqnarray}
which rules out of this case.

In the cases of (b1), (b2), (b3), we have $|k(\widetilde{\eta}_5-\widetilde{\eta}_6)+\widetilde{\tau}_5-\widetilde{\tau}_6|>k/4$ and then obtain the same estimate (\ref{count}), which also rules out of these cases.
%
%

Let us consider the case of (ii-2).
Since by $\max_{1\le m\le 5}m(\xi_m)\ll m(\xi_6)$,  we have $\widetilde{\eta}_m\in\{0,1,3,4\}, ~|\widetilde{\tau}_m|\ll k$ for $1\le m\le 5$, and $m(\xi_6)\lesssim k^{s-1/2}$.
We observe the case when $|\phi(\xi_1,\xi_2,\xi_3,\xi_4,\xi_5,\xi_6)|\gtrsim \langle k\rangle^2$.
An argument similar to above shows that this contribution to (\ref{rem-4}) is bounded by the sum of following three terms:
\begin{eqnarray*}
\frac{\langle k\rangle^{s-5/2}}{L^5}\left|\sum_{\xi_6\not\in{\cal A}_r}\sum^*a_{\xi_1}(t)\overline{a_{\xi_2}(t)}a_{\xi_3}(t)\overline{a_{\xi_5}(t)} a_{\xi_5}(t) m(\xi_6)\overline{e_{\xi_6}(t)}\right|,
\end{eqnarray*}
\begin{eqnarray*}
\frac{\langle k\rangle^{s-5/2}}{L^5}\int_0^t\sum_{\xi_6\not\in{\cal A}_r}\sum^*\left|\partial_{t'}\left(a_{\xi_1}(t')\overline{a_{\xi_2}(t')}a_{\xi_3}(t')\overline{a_{\xi_5}(t')} a_{\xi_5}(t')\right)m(\xi_6)\overline{e_{\xi_6}(t')}\right|\,dt'
\end{eqnarray*}
and
\begin{eqnarray*}
\frac{\langle k\rangle^{2s-3}}{L^5}\int_0^t\sum_{\xi_6\not\in{\cal A}_r}\sum^*\left|a_{\xi_1}(t')\overline{a_{\xi_2}(t')}a_{\xi_3}(t')\overline{a_{\xi_5}(t')} a_{\xi_5}(t')\partial_{t'}\overline{e_{\xi_6}(t')}\,dt'\right|,
\end{eqnarray*}
which are bounded by
\begin{eqnarray}
& & \frac{c}{\langle k\rangle^{5/2-s}}\nu^{5/2}\|e(t)\|_s+\frac{c}{\langle k\rangle^{5/2-s}}\nu^{9/2}\int_0^t\|e(t')\|_s\,dt'+\frac{c}{\langle k\rangle^{3-2s}}\nu^5 t\nonumber\\
& \lesssim &\frac{\nu^3}{\langle k\rangle^{5/2-s}}+    \frac{\nu^5}{\langle k\rangle^{3/2-s}}t+\nu^2\int_0^t\|e(t')\|_s^2\,dt'.\label{rem4411}
\end{eqnarray}
So we may suppose $|\phi(\xi_1,\xi_2,\xi_3,\xi_4,\xi_5,\xi_6)|\ll \langle k\rangle^2$.
By assumption of this case,
\begin{eqnarray*}
k(\widetilde{\eta}_1-\widetilde{\eta}_2+\widetilde{\eta}_3-\widetilde{\eta}_4+\widetilde{\eta}_5-\widetilde{\eta}_6)=o(k)+\widetilde{\tau}_6\in(-3k/4,3k/4)
\end{eqnarray*}
and
\begin{eqnarray}\label{sque}
\frac{1}{(2\pi)^2}\phi(\xi_1,\xi_2,\xi_3,\xi_4,\xi_5,\xi_6)=k^2\left(\widetilde{\eta}_1^2-\widetilde{\eta}_2^2+\widetilde{\eta}_3^2-\widetilde{\xi}_4^2+\widetilde{\eta}_5^2-\widetilde{\eta}_6^2\right)+O(k)
\end{eqnarray}
imply
$$
\widetilde{\eta}_1-\widetilde{\eta}_2+\widetilde{\eta}_3-\widetilde{\eta}_4+\widetilde{\eta}_5-\widetilde{\eta}_6=0,
$$
$$
\widetilde{\eta}_1^2-\widetilde{\eta}_2^2+\widetilde{\eta}_3^2-\widetilde{\eta}_4^2+\widetilde{\eta}_5^2-\widetilde{\eta}_6^2=0.
$$
Then Corollary \ref{s-1} allows us to obtain $\{\widetilde{\eta}_1,\widetilde{\eta}_3,\widetilde{\eta}_5\}=\{\widetilde{\eta}_2,\widetilde{\eta}_4,\widetilde{\eta}_6\}$, and any two of $\widetilde{\eta}_1,\widetilde{\eta}_3,\widetilde{\eta}_5$ do not coincide with each other.
Moreover, we have
$$
|\widetilde{\tau}_1-\widetilde{\tau}_2+\widetilde{\tau}_3-\widetilde{\tau}_4+\widetilde{\tau}_5-\widetilde{\tau}_6|\lesssim 1,
$$
which yields
$$
\langle \widetilde{\tau}_6\rangle\lesssim \max_{1\le m\le 5} \langle \widetilde{\tau}_m\rangle,
$$
and then $\max_{1\le m\le 5}m(\xi_m)\gtrsim m(\xi_6)$.
However this is out of the case.

\underline{Case (iii)}:
Since at least two elements of $\xi_m\not\in{\cal A}_r~(1\le m\le 5)$ satisfy $m(\xi_m)\gtrsim k^{s-1/2}$,  it is easy to see that this contribution to (\ref{rem-4}) is bounded by
\begin{eqnarray*}
c\left(\frac{1}{\langle k\rangle^{2(s-1/2)-s}}+1\right)\nu^2\int_0^t\|e(t')\|_s^2\,dt'\lesssim \nu^2\int_0^t\|e(t')\|_s^2\,dt',
\end{eqnarray*}
provided $s>1$.

Therefore we have the bound of (\ref{R3}) by
\begin{eqnarray}
& & -2\Im\sum_{\xi\in 2\pi\mathbb{Z}/L}\int_0^t m(\xi)^2\widetilde{R}_{\xi}^3(t')\overline{e_{\xi}(t')} \,dt'\nonumber\\
& \lesssim & 
o(\|e(t)\|_s^2)+\frac{\nu^3}{\langle k\rangle^{5/2-s}}+\frac{\nu^5}{\langle k\rangle^{3/2-s}}t+\nu^5 t^2+\nu^2\int_0^t\|e(t')\|_s^2\,dt'.\label{Rc}
\end{eqnarray}

%
%

\subsection{Proof of Proposition \ref{prop:approx}}
\indent

We now turn to the proof of Proposition \ref{prop:approx}.
The proof follows from continuity of solution with respect to $\|\cdot\|_s$-norm, together with the Gronwall inequality.

By Lemmas \ref{lem:men}, \ref{lem:R1}, \ref{lem:R2}, \ref{lem:R4}, and the estimate (\ref{Rc}), we get
\begin{eqnarray}\label{apr1}
\|e(t)\|_s^2\le A+Bt+Ct^2+c\nu^2\int_0^t\|e(t')\|_s^2\,dt',
\end{eqnarray}
where
$$
A=\frac{c\nu^3}{\langle k\rangle^{5/2-s}},\quad
B=\frac{c\nu^5}{\langle k\rangle^{3/2-s}},\quad
C=c\nu^5.
$$

We need to prepare the Gronwall inequality.
It is easy to prove that if the continuous function $F(t)$ satisfies
$$
F(t)\le A+Bt+Ct^2+c\nu^2\int_0^t F(t')\,dt',
$$
then we obtain that
$$
c\nu^2\int_0^tF(t')\,dt'\le Ae^{c\nu^2 t}-(A+Bt+Ct^2)+\frac{B}{c\nu^2}(e^{c\nu^2 t}-1)+\frac{2C}{c^2\nu^4}(e^{c\nu^2 t}-1-c\nu^2t),
$$
and then
$$
F(t)\le Ae^{c\nu^2 t}+\frac{B}{c\nu^2}(e^{c\nu^2 t}-1)+\frac{2C}{c^2\nu^4}(e^{c\nu^2t}-1-c\nu^2t).
$$
Applying the Gronwall inequality to (\ref{apr1}), we get
\begin{eqnarray*}
\|e(t)\|_s^2\lesssim \frac{\nu^3}{\langle k\rangle^{5/2-s}}e^{c\nu^2 t}+\frac{\nu^3}{\langle k\rangle^{3/2-s}}\left(e^{c\nu^2 t}-1\right)+\nu\left(e^{c\nu^2t}-1-c\nu^2t\right).
\end{eqnarray*}
This proves Proposition \ref{prop:approx}.
\qed

%
%

\section{Proof of Theorem \ref{thm-main}}\label{sec:thm}
\noindent

In this section we prove Theorem \ref{thm-main}.

\noindent
{\it Proof of Theorem \ref{thm-main}}.
We use a bootstrap (continuity) argument in a similar way to the proof of Lemma \ref{non-deg}.
By Lemma \ref{non-deg}, (\ref{K-fun}), (\ref{choice-K}), Lemma \ref{lem:pert} and Proposition \ref{prop:approx}, we in fact have $\|e(t)\|_s^2\lesssim \nu$ for $|t|\ll 1/(L^3\nu^2)$, and iterate Proposition \ref{prop:approx} to obtain Theorem \ref{thm-main}.
\qed

\end{document}